\documentclass[11pt, reqno]{amsart}

\usepackage{geometry,hyperref,amsrefs, enumitem, cancel, extarrows, color,lineno}   
\usepackage{changes}
\usepackage[mathscr]{euscript}

\numberwithin{equation}{section}

\theoremstyle{plain}
\newtheorem{theo}{Theorem}[section]
\newtheorem{lem}[theo]{Lemma}
\newtheorem{prop}[theo]{Proposition}

\newtheorem{cor}[theo]{Corollary}

\newtheorem{lemma}[theo]{Lemma}

\theoremstyle{definition}
\newtheorem{rem}[theo]{Remark}

\newtheorem{example}[theo]{Example}
\newtheorem{definition}[theo]{Definition}
\newtheorem{defi}[theo]{Definition}

\newenvironment{pf}{\noindent{\it Proof. }}{$\hfill\square$\par\medskip}

\theoremstyle{plain}

\theoremstyle{definition}

\newcommand{\beq}{\begin{equation}}
\newcommand{\eeq}{\end{equation}}
\newcommand{\beqn}{\begin{equation*}}
\newcommand{\eeqn}{\end{equation*}}
\renewcommand{\a}{\alpha}
\renewcommand{\b}{\beta}

\renewcommand{\d}{\delta}
\newcommand{\e}{\epsilon}

\newcommand{\f}{\varphi}
\newcommand{\g}{\gamma}
\newcommand{\h}{\eta}

\renewcommand{\k}{\kappa}
\newcommand{\lm}{\lambda}
\renewcommand{\o}{\omega}
\renewcommand{\q}{\vartheta}

\newcommand{\s}{\sigma}
\renewcommand{\t}{\tau}

\newcommand{\D}{\Delta}

\newcommand{\G}{\Gamma}

\renewcommand{\L}{\Lambda}
\renewcommand{\O}{\Omega}

\newcommand{\bC}{\mathbb{C}}
\newcommand{\bR}{\mathbb{R}}

\newcommand{\gF}{\mathfrak{F}}

\newcommand{\gX}{\mathfrak{X}}

\newcommand{\cD}{\mathscr{D}}
\newcommand{\cE}{\mathscr{E}}
\newcommand{\cF}{\mathscr{F}}

\newcommand{\cH}{\mathscr{H}}

\newcommand{\cK}{\mathscr{K}}
\newcommand{\cL}{\mathscr{L}}

\newcommand{\cU}{\mathscr{U}}

\newcommand{\cW}{\mathscr{W}}
\newcommand{\cX}{\mathscr{X}}

\newcommand{\p}{\partial}

\renewcommand{\square}{\kern1pt\vbox
{\hrule height 0.6pt\hbox{\vrule width 0.6pt\hskip 3pt
\vbox{\vskip 6pt}\hskip 3pt\vrule width 0.6pt}\hrule height0.6pt}\kern1pt}

\DeclareMathOperator\End{End}

\DeclareMathOperator\Id{Id}

\newcommand\Ker{\operatorname{Ker}}
\renewcommand\={:=}

\newcommand{\wt}{\widetilde}
\newcommand{\wh}{\widehat}

\newcommand{\bt}{\begin{theo}}
\newcommand{\et}{\end{theo}}
\newcommand{\bp}{\begin{prop}}
\newcommand{\ep}{\end{prop}}
\newcommand{\bc}{\begin{cor}\ \ }
\newcommand{\ec}{\end{cor}}
\newcommand{\bl}{\begin{lem}\ \ }
\newcommand{\el}{\end{lem}}
\newcommand{\bd}{\begin{definition}}
\newcommand{\ed}{\end{definition}}
\newcommand{\n}{\nabla}

\newcommand{\be}{\begin{equation}}
\newcommand{\ee}{\end{equation}}

\def\<#1,#2>{\langle\,#1,\,#2\,\rangle}
\newcommand{\arr}{\begin{array}{rlll}}
\newcommand{\ea}{\end{array}}
\newcommand{\bea}{\begin{eqnarray}}
\newcommand{\eea}{\end{eqnarray}}
\newcommand{\bean}{\begin{eqnarray*}}
\newcommand{\eean}{\end{eqnarray*}}

\newcommand{\Ga}[3]{\G_{#1 #2}^{\phantom{#1} #3}}
\newcommand{\GGa}[3]{{\mathbf \G}_{#1 #2}^{\phantom{#1} #3}}

\newcommand{\RR}[4]{\operatorname{R}_{#1 #2 #3}^{\phantom{#1 #2 #3} #4}}

\font\smallsmc = cmcsc9
\font\smalltt = cmtt8
\font\smallit = cmti8

\def\sideremark#1{\ifvmode\leavevmode\fi\vadjust{
\vbox to0pt{\hbox to 0pt{\hskip\hsize\hskip1em
\vbox{\hsize2cm\tiny\raggedright\pretolerance10000
\noindent #1\hfill}\hss}\vbox to8pt{\vfil}\vss}}}

\renewcommand{\sf}{shearfree }

\renewcommand{\e}{\operatorname{e}}

\newcommand{\ps}{{\operatorname{p}}}
\newcommand{\qs}{{\operatorname{q}}}

\newcommand{\grad}{{\operatorname{grad}}}
\newcommand{\Ric}{{\operatorname{Ric}}}
\newcommand\bex{\begin{example}}
\newcommand\eex{\end{example}}
\newcommand\br{\begin{rem}}
\newcommand\er{\end{rem}}

\makeatletter
\@namedef{subjclassname@2020}{%
  \textup{2020} Mathematics Subject Classification}
\makeatother

\title[Lorentzian manifolds with shearfree congruences]{Lorentzian manifolds with shearfree congruences\\ and K\"ahler-Sasaki geometry}
\author{Dmitri Alekseevsky, Masoud Ganji, Gerd Schmalz, Andrea Spiro}

\thanks{{\it Acknowledgments}. D.A. and G.S. were partially supported by the Australian Research Council, Discovery Grant DP130103485. D.A. was   supported  also  by  the grant  no. 18-00496S of
the Czech Science Foundation.}
\begin{document}
\definechangesauthor[color=red,name=Masoud]{M}

\maketitle
{\centering \ \\[-0.7cm] \it Dedicated to the memory of Alexandre Mikhailovich Vinogradov\\}
\begin{abstract}  We   study     Lorentzian manifolds  $(M, g)$  of  dimension  $n\geq 4$, equipped  with  a     maximally twisting shearfree null vector  field $\ps$, for which  the  leaf space  $S = M/\{\exp t\ps\}$  is a  smooth  manifold.   If  $n = 2k$, the quotient $S = M/\{\exp t\ps\}$  is naturally equipped with 
a subconformal  structure of contact type and,  in  the  most interesting  cases,  it is  a regular Sasaki manifold  projecting onto a  quantisable K\"ahler manifold of real dimension $2k -2$.  Going backwards through this line of ideas, for any  quantisable K\"ahler manifold   with  associated Sasaki manifold $S$,  we give the  local description of all   Lorentzian metrics $g$ on the   total spaces $M$ of   $A$-bundles  $\pi: M \to S$, $A = S^1, \bR$,  such that the generator  of the group  action  is a  maximally twisting shearfree $g$-null vector  field $\ps$.  We  also  prove  that on any  such Lorentzian manifold $(M, g)$  there exists   a non-trivial  generalized electromagnetic  plane wave  having $\ps$ as   propagating direction field, a result that can be  considered as a generalization of the classical $4$-dimensional Robinson Theorem.
We finally construct    a  2-parametric family of   Einstein metrics  on a trivial bundle  $M = \bR \times S$  for  any prescribed value of the Einstein constant. If $\dim M = 4$, the Ricci flat  metrics  obtained in this way  are  the well-known  Taub-NUT metrics.
\end{abstract}
\section{Introduction}
In this paper we   study  Lorentzian manifolds  $(M, g)$ of dimension $n \geq 4$,   equipped with  a  maximally twisting shearfree null vector  field $\ps$, for which  the  leaf space  $S = M/\{\exp t\ps\}$  is a smooth  manifold   and      $\pi : M \to S$  is a  principal bundle with one-dimensional structure group $A = \exp(t \ps)$ isomorphic  to  $ \bR$ or $S^1$.  In  case $n = 2k$, the  quotient $S = M/ A$  is odd-dimensional and naturally equipped with a  contact  distribution with a  subconformal  structure. In the  most interesting  situations, such subconformal structure  comes from a regular Sasaki structure on $S$ and the manifold $S$ is the total space of  an $A'$-bundle $\pi': S \to N$,  $A' = \bR, S^1$,  over a quantisable K\"ahler manifold $N$. Our main results concern  such Lorentzian manifolds fibering over regular Sasaki manifolds. More precisely, for any given quantisable K\"ahler manifold $N$ with  associated Sasaki bundle $\pi': S \to N$, we give the  local description of  the  Lorentzian metrics $g$ on the principle $A$-bundles  $\pi: M \to S$ over $S$,   for which the generator $\ps$ of the $A$-action  is a  maximally twisting shearfree $g$-null vector  field. We also    prove  a generalization of the classical $4$-dimensional Robinson Theorem,  namely  for any such Lorentzian manifold we prove the  existence of a non-trivial  generalized electromagnetic  plane wave, having  $\ps$  as the propagating direction field of the wave. We finally  construct   a  2-parametric family of   Einstein metrics  on any trivial bundle of the form $M = \bR \times S$  for any prescribed  value of the Einstein constant.  If  $\dim M = 4$,  the Ricci flat  metrics obtained in this way  are  the well-known  Taub-NUT metrics. 
\par
\smallskip
Let   $(M, g)$ be a $4$-dimensional Lorentzian space-time and $F$ a 2-form representing an  electromagnetic  plane wave, i. e. an harmonic decomposable   $2$-form  $F = \q \wedge \operatorname{e}^*$, determined by a null $1$-form $\q$ and  a  $g$-orthogonal   space-like  $1$-form   $\operatorname{e}^*$.  Any such electromagnetic field  is  associated with a flag structure, namely the pair of nested  distributions $\cK \subset \cW$,   given by  the spaces $\cK_x= \ker F_x \cap \ker ({\ast} F)_x$ and $\cW_x = \ker \q_x$, respectively. These distribution have the following physical interpretations:  $\cK$   is the  null  $1$-dimensional distribution  giving the propagating    directions   of the  wave and  $\cW$ is the   codimension  one distribution generated  by the wave fronts  and the propagating directions. Note that  a nested pair   $\cK \subset \cW$, with $\cK$  null and  $1$-dimensional  and $\cW$ of codimension one,  is the flag structure of an electromagnetic plane wav
 e 
only if  the  integral lines of  $\cK$  constitute a    {\it geodesic \sf congruence}, that is a family of curves $\g(t)$ tangent to $\cK$ and    such that   
$$\cL_{\dot \g} \cW \subset \cW\hskip 2 cm \text{and} \hskip 2cm  \cL_{	\dot \g} g = f g +\q\vee \h \quad\text{with}\quad  \q \= g(\dot \g, \cdot)$$  
for some function $f$ and a $1$-form $\h$. 
 This fact has a famous   converse,  the   Robinson Theorem (\cites{Ro, HM}):   {\it any geodesic \sf congruence  of a real analytic  Lorentzian $4$-manifold  $(M, g)$  locally coincides with the family  of propagation lines of a non-trivial electromagnetic plane wave}.  \par
 \smallskip 
An analogue of Robinson Theorem holds  for a large class of  Lorentzian manifolds $(M, g)$ of  higher even dimension $ 2k > 4$, provided that the classical notion of  electromagnetic plane wave is extended  as in the following definition, inspired by Trautman's discussion in   \cite{Tr2}. 
 A {\it generalised
electromagnetic plane wave} on $(M, k)$ is a 
 harmonic     $k$-form  $F = \q \wedge\a$, which is the wedge product of a null $1$-form $\q$ and a $(k-1)$-form $\a$ with the property that  the null vectors in $\cW = \ker \q$ are also in $\ker \a$.  Any such generalised plane wave determines a flag structure, namely the pair  $(\cK\= \ker F \cap \ker ({\ast} F), \cW\= \ker \q)$.   As in the $4$-dimensional case, a nested pair $\cK \subset \cW$ of  distributions might occur as a flag structure of a generalised electromagnetic plane wave only if it satisfies certain conditions. These constraints  are  satisfied if the distributions $\cK$ and $\cW$ are determined by some  geodesic \sf congruence on $(M, g)$.  Furthermore, an analogue of Robinson Theorem  for  generalised electromagnetic plane waves holds on any   even dimensional  Lorentzian manifolds of K\"aher-Sasaki  type, a very  large class of  space-times  which we present below.  \par
 \medskip
 These    facts, together with   their  tight relations  with   sub-Riemannian,  CR and  K\"ahler geometries discussed further in this paper, motivate  our interest for the  Lorentzian manifolds   equipped with  (geodesic) \sf congruences. We call them  {\it \sf Lorentzian manifolds}.  Such manifolds are particularly relevant  also because they provide  a large family of  examples of {\it null $G$-structures}, a class of structures on space-times  that has  recently received attention in the context of string and $M$-theory (see \cite{Pa} and references therein).  Note also that the  properties of  electromagnetic plane waves and the Robinson Theorem in higher dimension has recently  been object of intensive  investigations  (see e.g. \cites{Da, DPPR, Or, OPP, OPP1, OPP2, OPZ, OPZ1, SP}). \par
\smallskip
In  this paper we focus on a particular class of \sf Lorentzian   manifolds, the {\it regular} ones.  They are manifolds  on which the    \sf congruence   consists of   the  orbits of a one-dimen\-sional group $A = \bR$ or $S^1$ acting freely and properly.  For the manifolds of this kind,  the orbit space  $S = M/A$ is  a smooth manifold, which is  naturally  equipped with a subconformal structure $(\cD, [h])$, i.e. a pair given  by a codimension one distribution $\cD \subset TS$ and a positive definite conformal  metric $[h]$ on  such distribution. If in addition   $M$ is even dimensional and the \sf congruence satisfies  the so-called   {\it twisting condition},  the distribution  $\cD$  is   contact  and the subconformal structure is canonically associated  with a field of complex structures $J_x: \cD_x \to \cD_x$, $x \in M$,   that  makes the triple    $(S, \cD, J)$ a  strongly pseudoconvex almost CR manifold. \par
\smallskip
 The geometry of the orbit space $S = M/A$ of the regular \sf Lorentzian manifolds is  reacher and more interesting in case the  CR structure $(\cD, J)$ on $S$ is integrable (that is, with identically vanishing  Nijenhuis tensor) and  there is  a free proper action of a  one-dimensional group $A'$ of diffeomorphisms  preserving $(\cD, J)$ and a contact form $\theta$.  In fact, if these additional properties hold,  the almost CR manifold $(S = M/A, \cD, J)$ is a {\it regular Sasaki manifold}  and  the quotient $N = S/A' = M/(A \cdot A')$   is naturally equipped with  a  K\"ahler metric.   The regular \sf Lorentzian manifolds of this kind are called  {\it of K\"ahler-Sasaki type}. This is   the type of  space-times mentioned above,  admitting several  non-trivial  generalised electromagnetic plane  waves propagating along the lines of the \sf congruence. \par
 \smallskip
 The first  purpose of this paper is to discuss in  detail  each  of  the above mentioned   relations between    regular \sf Lorentzian manifolds,   strongly pseudoconvex CR manifolds,    Sasaki  manifolds  and  quantisable  K\"ahler manifolds (i.e.   K\"ahler manifolds that can be obtained as quotients of  some regular Sasaki manifold). For instance, we  fully characterise not only  the  subconformal structure $(\cD, [h])$  and the almost CR structure $(\cD, J)$   associated with a  given regular  \sf manifold, but, conversely,  also the  regular \sf structures  with  a prescribed subconformal structure or almost CR structure   on the obit space.  As a by-product,  we establish an exact procedure for   locally reconstructing all of  the  regular   \sf Lorentzian metrics  of K\"ahler-Sasaki type    projecting onto a  prescribed quantisable K\"ahler manifold and we obtain the  above mentioned higher dimensional version of the Robinson Theorem.  \par
\smallskip
The second  aim is   to show  that   the  proposed  construction
 of  regular  \sf Lorentzian manifolds  can be   used to describe  interesting families of  manifolds as,  for instance, new   classes  of  Lorentzian Einstein manifolds. \par
We prove that, for any choice of a real constant $\L$ and  of a  quantisable  $(n-2)$-dimensional   K\"ahler-Einstein manifold $N$ with  associated Sasaki manifold $S$,  there exists  a two-parameter family of Einstein  metrics  with Einstein constant $\L$  that make   $M = S \times \bR$ a regular \sf Lorentzian manifold of K\"ahler-Sasaki type.  These metrics are  explicitly  given in terms of the K\"ahler metric of $N$, the contact  form of the associated Sasaki manifold $S$ and two  functions of the fiber coordinate $t \in \bR$ of the trivial bundle $M = S \times \bR \to S$. These two  functions  are uniquely  determined  by the prescribed Einstein constant $\L$ of the   metric and two  arbitrary  real  constants $B$,  $C$ with $C > 0$.  We call such metrics {\it of Taub-NUT type}   since,  for the case    $N = \bC P^1 = S^2$, $S = S^3$ and $\L = 0$, the  corresponding   Ricci flat   metrics    are precisely  the $4$-dimensional Taub-NUT metrics on $S^3 \times \bR = \bR^4 
 \setminus \{0\}$.\par
We recall that  any  oriented  $2$-dimensional Riemannian manifold  $N$ is   K\"ahler and that,   if such a surface is   compact and with   integer  volume  form, then it  is  also a  quantisable  K\"ahler manifold.  Therefore, {\it for any prescribed real number $\L$,  our results  associate  a two-parameter  family  of Lorentzian Einstein manifolds $(M = \bR \times S, g)$ with Einstein constant $\L$,  admitting  electromagnetic plane waves,     to any given compact Riemann surface $N$ with integer volume form}.\par
\medskip
Shearfree congruences of null geodesics on 4-dimensional  space-times  have been  studied in the    physics literature  for a long time. In particular,  their relation
with the  strongly pseudoconvex   CR structure on the $3$-dimensional orbit spaces is well known and has been discussed in various contexts and
with diverse approaches. For an interesting and  stimulating overview, we refer to \cites{HLN, Tr1} and to  the extensive references therein.   Here  we tried to  give    a unified  and,  as much as possible,  exhaustive approach to  all of the cited relations between   such   four  types of important geometric structures,   the  shearfree structures,   the CR structures, the Sasaki  structures   and  the K\"ahler structures. We discuss them   in full generality, in arbitrary dimension,    using a    coordinate-free language. \par
\medskip
As we mentioned above,  the geometry of a  \sf Lorentzian  manifold can be considered as the geometry of a space-time,  in which there exists at least  one   electromagnetic plane wave.  However,  as  just a quick look at  the  starry sky    tells us,  the  Universe  is  pervaded by electromagnetic plane waves. This  fact together with the   Copernican Principle on the absence of privileged  points of the Universe suggests that any Lorentzian manifold    describing   a realistic  cosmological model  should  satisfy  the following
\begin{itemize}[leftmargin = 8pt]
\item[]
{\it  {\bf Copernican Principle of visual connectedness.} For any two points $x, y \in M$  with  $y$  in    the null  cone of $x$,  there is a  geodesic joining $x$ and $y$ which belongs to   a null shearfree congruence.}
 \end{itemize}
Maybe  a more realistic   conjecture  is  that such Copernican Principle is valid only locally,    for  sufficiently closed  points  $x$,$y$. Nonetheless we think that a local version of such principle  is   realistic and physically relevant. 
This also motivates an  interesting differential geometric problem:  characterise and possibly classify   the Lorentzian manifolds   satisfying  either  the above  principle or one of its local  variants.  We believe that the results of this paper  can be efficiently used  to attack this problem. We plan  to address  it  in a  future work.  \par
\medskip
The paper is structured as follows.  In  \S \ref{section1},  we introduce  and study  the notion of {\it \sf structure},  a geometric object  which  underlies  the  metric of a    \sf Lorentzian manifold. We then show how to reconstruct  all  Lorentzian metrics that are compatible with a given  \sf structure.   Note  that the classes of Lorentzian metrics compatible with a fixed   underlying \sf structure  are in natural one-to-one correspondence with the   {\it \sf optical geometries}  of  Robinson and Trautman  (\cites{RT, RT1}).  In \S \ref{section2}, we introduce the notion of   {\it regular}  \sf structures and we prove the one-to-one correspondence between these structures and the subconformal structures on the  orbit spaces. In \S \ref{section3} we  study  the   {\it twisting} regular \sf structures and  prove that, on a  fixed Lorentzian manifold $M$,  there is  a  bijection between this kind of  \sf structures and   the strongly pseudoconvex almost CR structures on the orbit space equipped with some  conformal class of positive endomorphisms.
At the end  of this section  we  introduce the regular   \sf structure of  K\"ahler-Sasaki type, we study their relations with  the  quantisable K\"ahler manifolds and we prove the advertised generalisation of Robinson Theorem.  In  \S \ref{section4}, we use the previous results
to determine the    Einstein Lorentzian metrics of K\"ahler-Sasaki type associated with a given K\"ahler-Einstein manifold and satisfying an appropriate ansatz.  In this way  we obtain the new family of  Lorentzian Einstein metric  mentioned above. In an  appendix, we give the explicit expressions for  the Christoffel symbols of a  metric of K\"ahler-Sasaki type associated with a prescribed quantisable    K\"ahler manifold.  Such expressions are
the outcomes of some   tedious but  very straightforward   computations and   are essential ingredients for  the  construction of the new metrics in   \S \ref{section4} and possibly of other types  of examples. \par
\medskip
\noindent {\it Acknowledgments.}  After completing our paper, we learned that  in \cites{FLT, TC}  A. Fino, T. Leistner  and A. Taghavi-Chabert,  
simultaneously and independently, obtained interesting results on higher dimensional 
shearfree congruences and Taub-NUT metrics, which partially overlap with and complement the results of this paper. We are 
sincerely grateful to Arman Taghavi-Chabert for bringing this to our attention. We also warmly  thank Marcello Ortaggio and  the anonymous referee for  useful observations and suggestions.\par
\medskip
\noindent{\it Notation.}
 The spaces of smooth real functions, vector fields  and  $1$-forms of  a $n$-dimensional manifold $M$ are denoted by $\gF(M)$,  $\cX(M)$ and $\O^1(M)$, respectively. Given  a tuple  of  vector fields $X_1, X_2, \ldots \in \gX(M)$, we indicate by $\langle X_1, X_2, \ldots \rangle \subset TM$ the distribution which they generate.  Given  a non-degenerate  metric $g$, for any   $X \in  \cX(M)$ we  denote   $X^{\flat} \= g(X, \cdot)$.
For any pair of $1$-forms  $\a, \b \in \O^1(M)$,   the symbol  $\a \vee \b$ stands for the symmetric tensor product $\a \vee \b  \= \frac{1}{2} \left(\a \otimes \b + \b \otimes \a\right)$.  
 We write $X \in \cW$ for any vector field $X$ tangent to a distribution $\cW \subset TM$ and  $\cL_X \cW \subset \cW$ if $X$ preserves $\cW$. 
\par
\medskip
\section{Shearfree  Lorentzian  manifolds and shearfree structures}\label{section1}
\subsection{First definitions}
 Let $\ps$ be a null vector  field  on  an  $n$-dimensional  Lorentzian manifold $(M,g)$. We associate  to $\ps$  the  following objects:
 \begin{itemize}[leftmargin = 15pt]
 \item[--] $ \cW:=\ps^{\perp}$ is the codimension  one  distribution  orthogonal to $\ps$;
 \item[--] the semipositive  degenerate metric $ h= g_\cW \= g|_\cW$   on $\cW$ which is  induced  by   $g$.
 \end{itemize}
 A null vector field $\ps$ is called {\it \sf} if  the (local) flow of $\ps$  preserves the subconformal structure
 $(\cW , [h])$, that is
\beq \label{twocond}      \cL_{\ps} X  \in   \cW \quad \text{for any} \ X \in\cW\ \quad \text{and}\quad   \cL_{\ps} h = f h\quad\text{for some}\ \   f \in \cF(M). \eeq
As  the next lemma shows,   the   condition  \eqref{twocond}  is equivalent  to
\beq \label{defsf2}   \cL_{\ps} g = f g + {\ps}^\flat \vee \h\qquad \text{for some function $f$ and a $1$-form}\ \ \h\ .\eeq
This   demonstrates  that  the  \sf vector fields can be considered as  generalisations of  the null conformal vector  fields. \par
\begin{lem}  \label{altdef}
  A   nowhere vanishing null  vector field  $\ps$ satisfies   \eqref{twocond} if and only if it satisfies \eqref{defsf2}.
  \end{lem}
\begin{pf} Assume that $\ps$ satisfies  \eqref{defsf2}.  Then  for any vector field $X \in \cW =   \ps^\perp$
$$  g(\ps, \cL_\ps X) = \ps (g(\ps,X) )    -  (\cL_\ps g)(\ps, X) - g([\ps, \ps], X)
= - f  g(\ps, X) - (\ps^\flat \vee \eta) (\ps, X) = 0\ ,
$$
showing that  $\cL_{\ps} X  \in  \cW $. From this  we also get
$$ \cL_\ps h  = (\cL_\ps g)_{\cW }
=  f g_{\cW  } +  (\ps^\flat \vee \eta)_{\cW  }  =
f h\ .$$
  Conversely, assume that  $\ps$ is a nowhere vanishing null vector field    satisfying   \eqref{twocond}. Then, around any point $x_o \in M$, we may consider a simply connected neighbourhood $\cU$ and a   vector field $\qs \in T\cU \setminus \cW|_{\cU}$
 so that  $g(\ps,\qs) =1 $. Let    $\cW' = \ker \qs^{\flat}_\cW$ be the kernel  of the $\qs^{\flat}$  in  $\cW|_{\cU}$.   This  determines the  following   direct sum decompositions of the tangent and the cotangent bundle of $\cU$
$$ T\cU = \langle \ps\rangle  + \cW' + \langle \qs\rangle  \ ,\qquad  T^*\cU = \langle \qs^{\flat}\rangle   + \cW'{}^* + \langle \ps^{\flat}\rangle\ .$$
These decompositions determine the following  direct sum decomposition of the bundle $S^2T^*\cU$ of the  symmetric  square
$$ S^2T^*\cU =  \langle \ps^\flat  \rangle \vee  \langle \ps^\flat \rangle  +  \langle \ps^\flat \rangle\vee \langle \qs^\flat \rangle + \langle \qs^\flat \rangle \vee  \langle \qs^\flat \rangle +  \langle \ps^{\flat} \rangle  \vee \cW'{}^* + \langle \qs^\flat\rangle \vee \cW'{}^* +\cW' \vee \cW'\ .$$
 Note that a symmetric $(0,2)$ tensor field  vanishes identically on any pair of vector fields in $\cW$ if and only if it takes values in $  \langle \ps^\flat  \rangle \vee  \langle \ps^\flat \rangle  +   \langle \ps^\flat \rangle\vee \langle \qs^\flat \rangle + \langle \ps^{\flat} \rangle  \vee \cW'{}^* $.
 On the other hand, by assumption, $(\cL_\ps g - f g)_{\cW|_{\cU}} = \cL_\ps h  - f h= 0 $. So, by the previous observation, at the points of $\cU$ we have that
 $\cL_\ps g - f g = \ps^\flat \vee \h$ for some (uniquely defined)  $1$-form $\h$.
 The uniqueness of $\h$ on   $\cU$  implies that \eqref{defsf2}  holds for   a unique  $1$-form $\h$ on $M$.
  \end{pf}
The   conditions  \eqref{twocond} are  mostly motivated by the fact that they correspond to the  two main properties of   the  propagating direction field $[\ps]$  of an electromagnetic plane wave in  General Relativity.  Indeed, on a $4$-dimensional space-time the first condition is equivalent to 
  the property that $\ps$ is geodesic, i.e.  $\n \ps = \lm \ps$ for a function $\lm$,   and  it encodes  the fact that  the photons 
   travel along null geodesics.  The second condition captures the fact that  the null field property of the electromagnetic plane waves (i.e., in terms of the electric and  magnetic fields,   $\overset{\rightarrow} E{\cdot} \overset{\rightarrow} H = 
 | \overset{\rightarrow}E|^2 - |\overset{\rightarrow}H|^2 = 0$) is preserved along their null propagation rays  (\cite{Ba}). It is also worth mentioning that if $\dim M = 3$ and  $\ps$ is a null vector field satisfying  just the first condition (i.e., it 
 preserves  the {\it $2$-dimensional} distribution $\cW = (\ps)^\perp \supset \langle \ps \rangle$),  then  the \sf   condition  $ \cL_{\ps} h= f h$, $h \= g_{\cW}$,  is   automatically satisfied.\par
\begin{definition} A manifold $M$  with a Lorentzian metric $g$ and  a \sf  vector field $\ps$ is called  {\it \sf Lorentzian manifold}.  The pair $(g,{\ps})$  is called  {\it \sf pair}.  The $1$-dimensional foliation of $M$,   which is determined
   by  the  integral curves of the $1$-dimensional distribution $\langle \ps \rangle$,   is called  {\it \sf congruence}.
\end{definition}
\par
\medskip
\subsection{Equivalent shearfree pairs. Standard and distinguished  pairs}
  \begin{definition}\hfill\par
  \begin{itemize}
 \item[i)]   Two   \sf pairs of the  form   $(g, \ps)$,  $(g' = \s g, \ps' = \t \ps)$  for some functions $\s > 0$,   $\t  \neq 0$,   are  called {\it equivalent}.   The   equivalence class of  $(g,\ps)$ is  denoted  by $([g], [\ps])$.
\item[ii)]   A \sf pair $(g,\ps)$  with autoparallel  vector  field  $\ps$,  i.e.  such that $\n_{\ps} \ps = 0$,  is called  {\it standard}.
 \item[iii)] A \sf pair  is called {\it  distinguished}  if
 \beq \label{ecconew} \cL_\ps g = \ps^\flat \vee \h\qquad \text{for some}\ \h \in \O^1(M)\ .\eeq
 \end{itemize}
 \end{definition}
\begin{prop}[\cites{RT, AGS}]\label{lemma23} \label{prop13} Let $(g, \ps)$ be a \sf pair. Then:
\begin{itemize}
\item[(i)]  $\ps $  is a geodesic vector field, i.e. $\nabla_\ps \ps= \lm \ps$ for some function  $\lm$. In other words, the corresponding \sf congruence is geodesic.
\item[(ii)]  Locally any \sf   pair  $(g, \ps)$    is  equivalent  to  a standard  pair  $(g, \ps') {=} (g, \t \ps)$. It is also locally equivalent  to   a  standard and  distinguished  pair  $(g', \ps') = (\s g, \t \ps)$. Such equivalent  standard and distinguished  pair is  uniquely determined  up  to   conformal factors $\s$, $\t$ that are constant along the $\ps$-orbits.
\end{itemize}
\end{prop}

\begin{pf}  (i) First of all, we  observe that   for any $g$-null vector field  $V$  the  following identity holds
\begin{multline} \label{a} g(\n_V V, X) = V(g(V, X)) -g(V, \n_V X) =  \\=( \cL_V g)(V, X) +  g(V, [V, X])   -g(V, \n_V X) =  (\cL_V g)(V, X) - g(V, \n_X V) =\\
=  (\cL_V g)(V, X) - \frac{X (g(V,  V))}{2}= ( \cL_V g)(V, X)\ .
\end{multline}
Since $\ps$ is null and shearfree,   this implies that for any vector field $X$
\beq\label{qua}  g(\n_\ps \ps, X) = \cL_\ps g(\ps, X)  = f g(\ps, X) + \frac{1}{2} \ps^\flat(X) \h(\ps) = (f +  \frac{1}{2}\h(\ps)) g(\ps, X) = \lambda g(\ps, X)\ . \eeq
From this  the  claim  (i) follows.  \par
\medskip
(ii)  Let
$(g' \= \s g, \ps' \= \t \ps)$ be a \sf pair  which is  equivalent to $(g, \ps)$  and
 denote by $\n$,  $\n'$ the Levi-Civita connections of $g$ and $g' = \s g$,  respectively.   Let also  $\ps^\flat \= g(\ps, \cdot)$ and $\ps'{}^\flat \= g'(\ps', \cdot)$. By Koszul's formula, for any   $X \in \gX(M)$
\begin{multline*} g'(\n_{\ps} \ps, X) =    \ps (g'(X, \ps) )-  g'([\ps, X], \ps)  =\\
= ( \ps(\s)) \ps^\flat(X) + \s \left(\ps(\ps^\flat(X) )-\ps^\flat([\ps, X])\right)
= ( \ps(\s)) \ps^\flat(X) + \s (\cL_{\ps}\ps^\flat)(X) \overset{\eqref{defsf2}}= \\
= \left( \ps(\s)
  + \s f +  \frac{1}{2}\s \h(\ps)\right) \ps^\flat(X)
\end{multline*}
for some function   $f $ and  a $1$-form $\h$.
Hence   $\n_{\ps} \ps =  \left( \ps(\s)
  + \s ( f +  \frac{1}{2}\h(\ps)) \right) \ps$ and
\beq
 \label{2.6}
  \nabla_{\ps'}\ps'  =  \t  \ps(\t) \ps + \t^2\nabla_{\ps}\ps      = \t  \bigg(  \ps(\t) + \t  (\ps(\s)  + \s    f)  +   \frac{1}{2}\s  \t  \h(\ps)  \bigg)\ps.
\eeq
On the other hand, $\cL_{\ps'} g =  \t \cL_{\ps}g - 2  \ps^\flat\vee d \t = \t f g + \ps^\flat \vee (\t \h - 2 d \t)$ and
\beq \label{2.8} \cL_{\ps'} g' =   \ps'(\s)g +  \s  \cL_{\ps'} g= \t \left( \ps(\s) +  \s  f\right) g   +  \ps^\flat\vee (\s \t \h - 2  \s d \t)\ .
\eeq
From \eqref{2.6} and  \eqref{2.8}, the local existence of a $\t$ and  a $\s$ so that $(g', \ps')$ is standard and/or distinguished  is a consequence of  the existence of  local solutions   to  the   system of differential equations
\beq \label{2.8bis}  \ps(\s) = -   \s f\ , \qquad  \ps(\t)  = -  \frac{1}{2} \s \t   \h(\ps) \ .  \eeq
This proves that any \sf pair is locally equivalent to a standard and distinguished pair. We leave to the reader the checking of the  last claim on uniqueness.
\end{pf}
\par
\medskip
\subsection{The shearfree  structure of a \sf  Lorentzian manifold}
We now introduce  a geometric object that   characterises the \sf  Lorentzian manifolds.    For this, we  observe  that
 an equivalence class  $([g], [\ps])$  of  \sf pairs determines the   codimension one distribution $\cW = \ker \ps^\flat$ and  the degenerate conformal metric $[h] = [g_\cW]$ on $\cW$,
that is the conformal  class of the semipositive degenerate metrics
induced by $g$.  The kernel  of such conformal metric is $\ker [h] = \langle \ps \rangle$ and the conformal metric  $[h]$ induces a positive definite  conformal metric on the so-called {\it screen bundle}  $\cW/\langle \ps \rangle$ (see \cite{RT}).  This motivates our
\begin{definition}  Let $M$ be an   $n$-dimensional  manifold  and  $(\cW, [h])$  a pair given by  a codimension one distribution $\cW \subset TM$ and  a conformal class $[h]$ of semipositive degenerate Riemannian  metrics on $\cW$ with  one dimensional kernel    $\cK_h = \ker h \subset \cW$.  \par
\begin{itemize}[leftmargin = 20pt]
\item[(i)]
The pair $(\cW, [h])$ is called  {\it \sf structure}  if it is preserved by   one (hence each)  vector field  $\ps$ in $\cK_h$, that is
\beq \label{2.9} \cL_\ps   \cW \subset \cW\qquad \text{and}\qquad  \cL_\ps [h] =  [h]\ .\eeq
Note that the pair $(\cW, [h])= (\ker \ps^\flat, [g_\cW])$ determined by an equivalence class $([g], [\ps])$ of   \sf pairs   is a   \sf structure.
\item[(ii)] A \sf pair $(g, \ps)$
 is called  {\it compatible with a  \sf structure $(\cW, [h])$}  if   $\ker \ps^\flat = \cW$ and  $[g_\cW]) = [h]$.  In this case  $g$ is called a {\it compatible metric} of the \sf structure.
 \end{itemize}
\end{definition}
In other words,   we may say that  a \sf structure  is a  semipositive degenerate  conformal metric  $[h]$ on  a codimension one  distribution $\cW$  satisfying  the following two conditions: (a) the  kernel $\cK_h = \ker h$ is one-dimensional and  (b)  $\cW$ is  invariant  under any vector field  $\ps \in \cK_h$.
\par
\medskip
 \subsection{Shearfree structures and optical geometries} \label{optgeo}
The notions of \sf structures and compatible metrics are tightly related with   the  optical geometries of Robinson and  Trautman (\cite{RT1}).    We recall that   an  {\it optical  geometry}  on a   manifold $M$ is a triple $(\cW, \cK, \{ g\})$ given by
 \begin{itemize}[leftmargin = 20pt]
 \item a codimension one distribution $\cW \subset TM$;
 \item a one-dimensional subdistribution $\cK \subset \cW$;
 \item an equivalence class $\{ g\}$ of Lorentzian metrics on $M$ satisfying the following  conditions:
 \begin{itemize}[leftmargin = 15 pt]
 \item  for any $g \in \{g\}$, one has  $g(\ps, \cW) = 0$ for any vector field $\ps \in \cK$;
 \item any two metrics $g, g' \in \{g\}$  are related   by
   $g' = \s g + \a \vee \b$  for some positive  function $\s : M \to \bR$   and  $1$-forms  $\a, \b$     with $\a \neq 0$ and  $\a|_{\cW} = 0$.
\end{itemize}
 \end{itemize}
 \par
 \smallskip
Now, any   \sf optical geometry  $(\cW, \cK, \{g\})$  (i.e. any optical geometry admitting   a \sf vector field $\ps \in \cK$) determines an associated \sf structure, namely the pair $(\cW, [h] = [g_\cW])$,  $g \in \{g\})$. In this case,   the  class  $ \{g\}$
 consists  of all of the  compatible metrics  of   such  \sf structure.
Conversely,  any \sf structure $(\cW, [h])$  determines   the  \sf  optical geometry $(\cW, \cK = \cK_h, \{g\})$ 
 with  class  $\{g\}$  given by  the  compatible metrics of the \sf structure.
\par
\medskip
\subsection{Reconstructions of  \sf optical geometries from  \sf structures}
\label{section15}
 We now  focus on the problem of reconstructing  a   \sf optical geometry $(\cW, \cK, \{g\})$ starting from the  associated \sf structure $(\cW, [h])$, that is of determining  all \sf metrics and \sf pairs  that are compatible with such \sf structure.\par
  \bigskip
     Given a codimension one distribution $\cW \subset TM$,   a $1$-form $\q$ such that $ \ker{\q} = \cW$ is called  a   {\it    defining $1$-form for $\cW$}.
   Note that,  for  any  vector field $\qs$ which is transversal  to $\cW$,     there is   a  unique defining $1$-form $\q$  such that  $\q(\qs) = 1$.
    \par

 \begin{definition} Let $(\cW, [h])$ be a \sf structure on a manifold $M$ and $\cK_h = \ker h \subset \cW$ the corresponding  one-dimensional kernel subdistribution.  A {\it  rigging} for $(\cW, [h])$ is a pair $(\cW', \qs)$ given by a subdistribution $\cW'\subset \cW$, which is complementary to  $\cK_h$ in $\cW$, and a vector field $\qs$ which  is transversal to $\cW$.
 \end{definition}
 Notice that  any pair given by  a  rigging  $(\cW', \qs)$   and a vector field $\ps \in \cK_h$ determines a direct sum decomposition
 $$    TM = \cK_h + \cW' + \langle \qs \rangle $$
 and two  $1$-forms $\ps^*$,  $\qs^*$ satisfying  the  conditions
 \beq \ps^*(\ps) = 1  \ ,\ \   \ker \ps^* = \cW' + \langle \qs\rangle\qquad \text{and}\qquad \qs^*(\qs) = 1  \ ,\ \  \ker \qs^* = \cW' + \langle \ps\rangle = \cW \ .\eeq
 Note also that $\qs^*$ coincides with  the  defining $1$-form $\q =  \qs^*$ for $\cW$   with  $\q(\qs) = 1$ and it   is  thus uniquely determined by  $\cW$ and  $\qs$.\par
 \smallskip
 With  a  rigging $(\cW', \qs)$ and a  degenerate metric  $h \in [h]$ we associate  the  unique  Lorentzian  metric $g$ that satisfies the
   conditions $g(\ps, \qs) = 1$,  $g_\cW = h$,   $g(\ps, \ps) = g(\qs, \qs) = 0$  and $\cW' = \langle \ps, \qs\rangle^\perp$.   This metric is
\beq\label{reconstructed}  g = h  + 2 {\qs}^*\vee {\ps}^*  =  h  + 2 \q \vee  {\ps}^* \ ,\eeq
where   $h$ is  considered  as a degenerate metric on $TM$  with  kernel  $\ker h =  \cK_h + \langle \qs\rangle$.   \par
These observations yield  the  following theorem, which  gives a complete description of all optical geometries  that are associated  with a  \sf structure.
\begin{theo} \label{compatibility}  Let  $(\cW, [h])$ be a \sf structure on $M$.
\begin{itemize}[leftmargin = 20 pt]
\item[(1)] For any triple $(h, \ps, (\cW', \qs))$  formed  by
\begin{itemize}[leftmargin = 25 pt]
 \item[(a)] a degenerate metric  $h \in [h]$,
  \item[(b)] a nowhere vanishing vector field $\ps \in \cK_h$  and 
  \item[(c)] a rigging  $(\cW', \qs)$,
  \end{itemize}
the corresponding metric   \eqref{reconstructed} together with the vector field $\ps$ gives a  \sf pair $(g, \ps)$ which is compatible with $(\cW, [h])$.
\item[(2)]Conversely, if $(g, \ps)$ is a  compatible \sf pair for $(\cW, [h])$, then locally there is a  rigging $(\cW', \qs)$
 such that $g$  is the metric  \eqref{reconstructed} determined by some triple  $(h = g_\cW, \ps, (\cW', \qs))$. Such a rigging can be determined in two different ways: 
 \begin{itemize}[leftmargin = 15pt]
 \item[(a)]   either   by  choosing   a  subdistribution $\cW' \subset \cW$  which is  complementary  to $\cK$  and  setting   $\qs$  to be the   vector field with       $g(\qs, \cW' + \langle \qs\rangle) = 0$ and $g(\qs, \ps) = 1$
  \item[(b)] or   by  choosing a   vector field $\qs$ which is transversal to $\cW$ and  setting $\cW' = \langle \ps, \qs\rangle^\perp$.
  \end{itemize}
\end{itemize}
\end{theo}
\begin{pf} (1) By construction, the vector field $\ps$ is a  nowhere vanishing null vector field of $g$ with  $\ker \ps^\flat = \cW$. Moreover, $g_\cW = h$ and $\ps \in \cK_h = \ker h$. Since $(\cW, [h])$ satisfies \eqref{2.9} for any vector field in $\cK_h$, it follows that $\ps$ is \sf for $g$. \par
    (2)  Let  $(g, \ps)$  be a compatible pair for  $(\cW, [h])$, i.e.  such that  $\ker \ps^\flat = \cW$,  $g_\cW \in [h]$ and $\ps \in \cK_h$. Then, for any    vector  field  $\qs$ such that $g(\ps, \qs) = 1$,  the corresponding  codimension two distribution $\cW' = \langle \ps, \qs\rangle^\perp$ defines a rigging $(\cW', \qs)$  for  $([h], \cW)$ and the  associated  metric  $g$ has the form \eqref{reconstructed}.  The same holds   for any choice of a  subdistribution $\cW' \subset \cW$,  which is complementary to $\langle \ps \rangle$,  and   the    vector field characterised  by  the conditions  $g(\qs, \cW' + \langle \qs\rangle) = 0$ and $g(\qs, \ps) = 1$. \end{pf}
The following corollary shows that the collection of all compatible metrics for  a given \sf structure are locally  parametrised by the pairs formed by a positive function $\s$ and a
$1$-form $\varpi$   which does not vanish on the distribution $\cK_h$. \par
\begin{cor}    Let   $(\cW, [h])$ be a \sf structure on $M$ and $(g_o, \ps_o)$  be  a  compatible \sf pair,   hence of the form
\beq\label{reconstructed-bis}  g_o = h_o  + 2  \qs_o^*\vee \ps_o^*  \ ,  \eeq
where $\qs_o^*$ and $\ps_o^*$  are the $1$-forms determined by a triple $(h_o, \ps_o, (\cW'_o, \qs_o))$  as in Theorem \ref{compatibility}, and  
 $h_o \in [h]$  is considered   as  a degenerate metric on  $TM$ with $\ker h_o = \langle \ps_o, \qs_o \rangle$.
Then any  other  compatible metric $g$  has locally the form
\beq \label{parame} g =  \s h_o  + 2 \qs^*_o\vee \varpi=   \s h_o  + 2 \ps^\flat_o\vee  \varpi\qquad \text{with}\ \ps_o^\flat = g_o(\ps_o, \cdot) =\qs_o^*\eeq
for some positive function   $s>0 $    and  a   $1$-form $\varpi$ with  $\varpi(\ps_o) \neq 0$.
\end{cor}
\begin{pf} By Theorem \ref{compatibility},
if $g$ is a compatible metric, then there is a triple   $(h , \ps_o, (\cW', \qs))$,  with $h = \s h_o \in [h]$,  such that  $g$  has the form
$g = \s \wt h_o  + \qs^*\vee \wt \ps_o^*$, where:
\begin{itemize}[leftmargin = 20pt]
\item[--]  $\wt h_o$  is the extension of $h_o$ as a degenerate metric on $TM$ with $\ker \wt h_o = \langle \ps_o, \qs\rangle$ and
\item[--]$\qs^*$ and $\wt \ps^*_o$   are the $1$-forms defined by
 $$ \wt \ps_o^*(\ps_o) = 1  \ ,\  \ker \ps_o^* = \cW' + \langle \qs\rangle\qquad \text{and}\qquad \qs^*(\qs) = 1  \ ,\  \ker \qs^* = \cW' + \langle \ps_o\rangle = \cW \ .$$
 \end{itemize}
It follows that
\begin{itemize}[leftmargin = 25pt]
\item[(a)] there exists a function $\lambda$ such that $\qs^*  = \lambda \qs_o^*$;
\item[(b)]  $\langle \ps_o\rangle \subset \ker (\wt h_o  - h_o)$ and hence   $\wt h_o - h_o = 2 \qs_o^* \vee \varpi_1$ for some $1$-form $\varpi_1$.
\end{itemize}
Therefore
$g = \s h_o  + 2 \qs_o^*\vee \varpi$ with $\varpi =  \varpi_1 + \lambda \wt \ps_o^*$ and \eqref{parame} holds.  \par
Conversely, if $g$ has the form \eqref{parame}, then it is a Lorentzian metric for which $\ps_o$ is a null vector field,  $\ps_o^\perp = \cW$  and $\ps_o$  is  a \sf vector field  preserving the  \sf structure  $(\cW, [h_o])$. This means that $g$ is a compatible metric.
  \end{pf}
We conclude this section   by giving  a local explicit description of the compatible metrics for  a given \sf structure $(\cW, [h])$ in terms of some   frame field  aligned with the distributions  $\cW$ and  $\cK_h \subset \cW$. \par
\smallskip
Let
$g_o$ be a  compatible metric  for the \sf structure  $(\cW, [h])$ and  let   $(h_o, \ps_o, (\cW', \qs_o))$   be a triple, formed  by a degenerate metric  $h_o \in [h]$,  a vector field $\ps_o \in \cK_h$ and a rigging  $(\cW'_o, \qs_o)$, which determines $g$  as   in (2) of Theorem \ref{compatibility}.  Consider also a  (local) frame field on $M$ of the form $(\ps_o, e_1, \ldots, e_{n-2}, \qs_o)$  where the  vector fields $e_i$ are in  the subdistribution $\cW'_o\subset \cW$. We  denote by $(\ps_o^* , e^1,\dots, e^{n-2}, \qs^*_o)$ the  dual coframe field.  Any other rigging $(\cW', \qs)$ for the   \sf structure  $(\cW, [h])$  has the form
\beq \label{52bisbis} \qs \=a \qs_o + b \ps_o + c^i  e_i\ ,\qquad \cW' = \langle e_i + d_i \ps_o\rangle\eeq
where $a, b, c^i, d_j$ are  smooth functions  with  $a \neq 0$ at all points.   However, since  by claim  (2a) of Theorem \ref{compatibility},  the  subdistribution  $\cW'$ can be fixed arbitrarily, {\it  with no loss of generality from now on we  assume that $\cW' = \cW'_o$ and    $d_j \equiv 0$}. \par
\medskip
The  $1$-forms $\wt \ps_o^*$ and $\qs^*$   satisfying the conditions
 \beq \wt \ps_o^*(\ps_o) = 1  \ ,\ \   \ker \ps_o^* = \cW'_o+ \langle \qs\rangle\qquad \text{and}\qquad \qs^*(\qs) = 1  \ ,\  \ker \qs^* = \cW'_o+ \langle \ps_o\rangle = \cW\ .\eeq
 can be  expressed in terms of the coframe field  $(\ps^*_o, e^i, \qs^*_o)$ as
\beq \wt \ps^*_o = \ps_o^* - \frac{b}{a} \qs_o^*\ ,\qquad \qs^* = \frac{1}{a} \qs^*_o\ .\eeq
However,  for any  $h = \s h_o \in [h]$, the corresponding   degenerate extension to $TM$   with $\ker h = \langle \ps_o, \qs \rangle$ is equal to
\beq h = \s  h_{ij} e^i \vee e^j - 2\frac{\s c^i h_{ij}}{a} \qs_o^* \vee e^j +   \frac{\s h_{ij} c^i c^j}{a^2} \qs_o^* \vee \qs_o^*\qquad h_{ij} \= h_o(e_i, e_j)\ .\eeq
 Since  any compatible  $g$ is as in \eqref{reconstructed}  for some  triples $(h = \s h_o, \ps_o, (\cW', \qs))$, we  conclude that
{\it  in terms of the  coframe field
$(\ps_o^*, e^i, \qs_o^*)$ any compatible metric  has the form}
 \begin{multline} \label{buona} g =  \s h_{ij} e^i\vee e^j + \qs^*_o \vee \left(\frac{2}{a} \ps^*_o  - 2 \frac{\s c^i  h_{ij} }{a} e^j +   \frac{\s h_{ij} c^i c^j}{a^2} \qs_o^* - \frac{2b}{a^2} \qs_o^*\right) = \\
 = \s\left( h_{ij} e^i\vee e^j + \qs^*_o \vee \left(\a\ps^*_o   + \g^i h_{ij} e^j  + \b  \qs_o^*\right) \right)
,\end{multline}
\beq \label{buona-bis}  \text{where}\ \a \= \frac{2}{a \s}\ ,\quad \g^i \= - \frac{2 c^i}{a}  \ ,\quad \b \= - \frac{2 b}{a^2\s}+  \frac{h_{ij} c^i c^j}{a^2}\ .\eeq
This gives a  parameterisation of all compatible metrics in terms of the $(n+1)$-tuple of functions $(\s, \a, \b, \g^i)$, which are  in turn  determined by arbitrary  functions $\s> 0, a \neq 0, b, c^i$,  as it is indicated in \eqref{buona-bis}.\par
\smallskip
Conversely,   any Lorentzian metric of the form  \eqref{buona}  for some  tuple $(\s, \a, \b, \g^i)$ with  $\s > 0$,  $\a \neq 0$,  is  a compatible  metric for $(\cW, [h])$, since it is associated with  the triple $(h = \s h_o, \ps_o, (\cW'_o, \qs))$,  in which  the vector field $\qs$ of the  rigging
$(\cW'_o, \qs)$ is defined by
\beq
\begin{split}  \qs &= a \qs_o + b \ps_o + c^i e_i\\
&\text{with}\ \  a \= \frac{2}{\a \s}\ ,\qquad  b \=  \frac{2}{\a \s} \left(- \b +  \frac{\g^\ell \g^m  h_{\ell m } }{4}  \right) \ ,\qquad   c^i \= - \frac{\g^i}{\a \s}\ .
\end{split}
\eeq
\par
\medskip
\subsection{Twisting \sf  structures}
Let $(\cW,[h])$ be a \sf structure on  $M$.  The {\it degree at  $x \in M$} is the integer
  $d_x(\cW):= \dim \Ker d\q|_{\cW_x}$   for some (hence,  for any) defining $1$-form $\q$ for $\cW$.   The \sf structure  is said to be  of {\it constant degree $d$} if the
  degree is  equal to $d$ at all points.  \par
  \begin{lem}  For any \sf structure $(\cW,[h])$ with  defining $1$-form $\q$, the one dimensional distribution $\cK_h$ is in $\ker d \q_{\cW}$. In particular,
 if  $\dim M = 2m$,  the minimal possible degree   is   $1$ and in this case  $\cK_h =\ker d \q|_{\cW}$. If  $\dim M = 2m +1$,  the minimal possible degree is $2$.
  \end{lem}
  \begin{pf} Since  ${\ps} \in  \cK_h$ preserves  $\cW$,  we have
 $ \lambda \theta = \cL_{\ps} \q  = d \q({\ps},\cdot)  + d(\q(\ps)) =$ $d \q({\ps},\cdot) $. Hence $d \q ({\ps}, \cW)= 0$.
\end{pf}
\bd  A \sf structure $(\cW,[h] )$  is called {\it twisting} if it  has     constant minimal  degree.
\ed
\par
\medskip
\section{Compatible \sf metrics of Robinson-Trautman bundles} \label{section2}
  In this  section we consider  our  main objects  of interest,   the {\it regular} \sf pairs $(g,\ps)$  and  the  associated   \sf structures, which we  call     {\it Robinson-Trautman   structures}.   Our aim is  to  characterise  these \sf structures  and  the  corresponding  optical geometries.   \par
\subsection{Regular \sf pairs and Robinson-Trautman structures}\hfill\par
  \bd A \sf pair $(g, \ps)$ on a manifold $M$  is called {\it regular} if the vector field $\ps$ is complete
  and generates a   1-parameter group   $A = e^{t{\ps}}$ of diffeomorphisms (isomorphic  to $\bR$ or $S^1$)  acting freely and properly  on  $M$.  \ed
  It is known that in this case the orbit space $S = M/A$ is a smooth manifold and  the   quotient map   $\pi : M \to S = M/A$  is  a smooth principal $A$-bundle.  The  corresponding \sf structure $(\cW= \ker \ps^\flat, [h] = [g_{\cW}] )$ is $A$-invariant   and  the  kernel distribution $\cK_h$ coincides  with the vertical distribution $T^{\operatorname{v}} M \subset TM$ of the bundle. This motivates    the following
  \bd
 A  \sf  structure  $(\cW, [h])$ on the total space of a principal  $A$-bundle $\pi: M \to S$, $A = \bR$ or $S^1$, is called  {\it Robinson-Trautman (RT) structure}
if it is $A$-invariant and  the  kernel distribution  is  $\cK_h = T^{\operatorname{v}} M$.  A principal $A$-bundle  $\pi: M \to S$ equipped   with an RT-structure $(\cW, [h])$ is called {\it Robinson-Trautman (RT) bundle}.
\ed
 We remark  that on an RT bundle $(\pi: M \to S, (\cW, [h]))$ the  fundamental vector  field   $\ps_o$ which  corresponds to the element $1 \in Lie(A) = \bR$, that is   the   velocity vector field  of  the $1$-dimensional Lie  transformation group  $A = \bR$ or  $S^1$,  gives   a canonical section  of $\cK_h = T^{\operatorname{v}}M$.
\par
\medskip
\subsection{Global standard  and distinguished   \sf pairs on RT bundles}
 The following proposition can be considered as a global version of the claim (ii) of Proposition \ref{lemma23} on the  RT bundles with structure group $\bR$.
\bp \label{proposition33}
On any  RT bundle   $(\pi: M \to S,(\cW,[h]))$   with    structure group $A = \bR$ 
  there exists  a  global   compatible   \sf pair $(g, \ps)$,  which  is standard  and distinguished,  i.e.  with
  $\nabla_\ps \ps =0$   and   $\cL_\ps g = \ps^{\flat} \vee \h$ for some globally defined $1$-form $\h$.
\ep
\begin{pf}
Let $\{\cU_\alpha\}$ be a trivialising cover of the base manifold $S$. Then for each $\cU_\a$ we have that
 $M|_{\cU_{\alpha}}\simeq \cU_{\alpha}\times \bR$ and we may consider  a \sf pair $(g_o, \ps_o)$ in which  $\ps_o=\partial_t$ for an appropriate  fiber coordinate $t$. By the proof of Proposition \ref{prop13},  if  the   $\cU_\a$ are sufficiently small, for  each of them
 there is  a pair $(\s_\a, \t_\a)$ of  strictly positive functions that are solutions
to the system  \eqref{2.8bis} over   $M|_{\cU_\a} = \cU_\a \times \bR$. Let $\{\chi_{\alpha} \}$ be a partition of
unity which is subordinated to  the open cover   $\{\cU_\alpha\}$ of $S$ and $\{\wt \chi_\a = \pi^* \chi_\a\}$  the corresponding family of  pulled-back functions on $M$. Notice that each   $\wt \chi_\a$ is constant along the integral curves of the vector field $\ps_\a \= \s_\a \ps_o$. This implies that  the pairs $(\wt \chi_\alpha {\cdot} \s_\a,  \wt \chi_\a {\cdot}\t_\a)$ are solutions to   \eqref{2.8bis}.  It follows  that
$(\s, \t) \= \left(\sum_\a \chi_\alpha \s_\a, \sum_\a \chi_\a \t_\a\right)$
is a global solution to the system \eqref{2.8bis} and determines a global standard and distinguished pair
$(g = \s g_o, p = \t p_o)$.\end{pf}
\begin{rem}  
In Proposition \ref{proposition33} the assumption   $A =  \bR$  is essential. In fact, by considering appropriate quotients of such an RT bundle, one can produce  examples of   RT bundles  with structure group $A = S^1$ admitting  no  global  standard or distinguished compatible \sf pair.
\end{rem}
\par
\medskip
\subsection{RT structures  and   subconformal structures}
\begin{definition}
A {\it sub-Riemannian} (resp.  {\it  subconformal}) {\it structure}  on  a manifold $S$  is  a  codimension one distribution $\cD \subset TS$ equipped   with  a Riemannian metric $g^{\cD}$  (resp.  a  conformal metric $[g^{\cD}]$) (\footnote{We  do not assume that $\cD$ is bracket generating,  i.e. a  contact distribution.  The  structures with  this additional  hypothesis are discussed  in  \S \ref{section331} below.}).
\end{definition}
The next proposition gives a fundamental relation  between the  RT structures and  the  subconformal structures.
\bp There is a natural one-to-one correspondence between the  RT structures $(\cW, [h])$ on the total space of  a  principal $A$-bundle  $\pi: M \to S$, $A = \bR$ or $S^1$,  and the  subconformal structures $(\cD, g^\cD)$ on   the base manifold $S$.
\ep
\begin{pf}    Since the codimension one distribution $\cW \subset M$ is  $A$-invariant, it projects onto  a codimension  one  distribution $ \cD:= \pi_* (\cW)$ on $S$. By a similar reason,    the   degenerate $A$-invariant conformal  metric $[h]$ on $\cW$  projects onto
a  conformal   metric $[g^{\cD}] = \pi_*([h])$ on $\cD$. This  associates  a   subconformal structure  $(\cD, [g^{\cD}])$  on  $S$ with any RT structure on $M$. Conversely,   given  a  subconformal structure  $(\cD, [g^\cD])$ on  $S$.  Then: 
\begin{itemize}[leftmargin = 20 pt]
\item[(a)] the preimage
$ \cW :=\pi_*^{-1}(\cD)$ is  an $A$-invariant codimension one  distribution  on $M$ which  contains $T^{\operatorname{v}}M$; 
\item[(b)]  the  pull-back   $[h] := (\pi^*)[ h^{\cD}]$  is  an $A$-invariant degenerate conformal metric  on $\cW$ with   kernel $\cK_h = T^{\operatorname{v}}M $.
 \end{itemize}
 In particular, $(\cW, [h])$ is an RT structure on $\pi: M \to S$.
\end{pf}
\par
\medskip
\subsection{$A$-invariant   metrics on RT bundles}
  Let  $\pi: M \to S$ be an RT bundle with structure group $A = \bR$ or $S^1$. Denote by  $(\cW, [h])$ and     $(\cD = \pi_*(\cW), [g^\cD] = \pi_*([h]))$  the corresponding RT structure and     subconformal structure on $M$ and $S$, respectively.
Let also  $\ps_o \in T^{\operatorname{v}} M$  be the   fundamental  vector of such principal bundle.  \par
  \smallskip
   Consider  a principal connection $\cH \subset TM$ on $M$,  that is  an  $A$-invariant (horizontal) distribution  complementary to the vertical distribution $T^{\operatorname{v}} M \subset TM$. We
 recall that:
 \begin{itemize}[leftmargin = 20 pt]
 \item[(a)]   any vector field $Y \in \gX(S)$ has a unique $A$-invariant horizontal lift $Y^{\operatorname{h}}$ in  $\cH \subset TM$   projecting onto $Y$;
 \item[(b)]  since the   kernel subdistribution $\cK_h \subset \cW$ is equal  to $\cK_h = T^{\operatorname{v}} M$, the intersection  $\cW' = \cW \cap \cH$ is an  $A$-invariant   subdistribution  complementary to  $\cK_h$.
 \end{itemize}
 It follows that   for any    vector field  $Z  \in \gX(S)$  that  is transversal to $\cD$, the pair $(\cW' {=} \cW \cap \cH$, $\qs {=} Z^{\operatorname{h}})$ is an $A$-invariant  rigging for the  RT structure $(\cW, [h])$.  Hence,
 given a  sub-Riemannian metric $g^\cD \in [g^\cD]$ on the distribution  $\cD \subset TS$, the  triple   $(h = \pi^*(g^\cD),\ps_o,$ $(\cW', \qs = Z^{\operatorname{h}}))$ is  $A$-invariant and determines an $A$-invariant compatible  metric $g$ by (1) of Theorem \ref{compatibility}. This  gives  the  following useful result. \par
   \bt  Let $\pi: M \to S $ be  an RT bundle with fundamental vector field $\ps_o$,    RT structure $(\cW, [h])$ on $M$ and    subconformal structure $(\cD, [g^\cD])$ on $S$  as above. Any  pair $(g^\cD, Z)$, formed by a sub-Riemannian metric $g^\cD \in [g^\cD]$ and a $\cD$-transversal vector field  $Z \in \gX(S)$,  determines  the  $A$-invariant compatible Lorentzian metric on $M$
    \beq \label{pupo}  g \= \pi^*(g^\cD) + \q  \vee \ps_o^* \qquad \text{where} \eeq
 \begin{itemize}[leftmargin = 18pt]
 \item[(a)]  the tensor $\pi^*(g^\cD)$ is considered as  a degenerate metric on $TM$ with kernel  $\langle \ps_o, Z^{\operatorname{h}}\rangle$,
 \item[(b)]  $\q$ is the   defining $1$-form  for the distribution  $\cW = \langle \ps_o \rangle  + \cW'$  such that    $\q(Z^{\operatorname{h}})=1$;
 \item[(c)] $\ps_o^*$ is   the $A$-invariant defining  form of    $\cH = \cW' + <\qs>$  such that
       $  \ps_o^*(\ps_o) = 1$.
  \end{itemize}
   The \sf pair $(g, \ps_o)$  is  standard  and distinguished, actually $\ps_o$ is Killing.
  \et
  \par
  \medskip
\section{Twisting Robinson-Trautman bundles of K\"ahler-Sasaki type} \label{section3}
Let $M$ be  the total space    of an even dimensional {\it twisting}  RT bundle $(\pi: M \to S,(\cW, [h]))$ with structure group $A = \bR, S^1$.  As we  will  see, the assumption that the \sf structure is twisting is equivalent to the hypothesis that  the corresponding  subconformal structure  $(\cD = \ker \vartheta , [h^{\cD}])$  on   $S$ is of contact type. This subconformal structure canonically determines   a strongly pseudo-convex almost CR  structure and, if  appropriate regularity conditions are satisfied, such CR structure makes $S$ 
 the total space of an   $A'$-bundle $\pi^S: S \to M$, $A' = \bR$ or $S^1$,  over a quantisable K\"ahler manifold $(N, J^N, g^N)$. \par
 In   this section   we   review  the main  definitions and certain  basic facts on subconformal structures  of contact type, Sasaki manifolds and quantisable K\"ahler manifolds. After that we  establish  in detail  the exact  relations   between the  twisting RT bundles  and  such  geometric structures.
\par
\medskip
\subsection{Subconformal structures of  contact type and CR  structures}
 \subsubsection{Sub-Riemannian  and  subconformal structures of  contact type} \label{section331}
 \begin{defi} Let $S$ be an odd dimensional manifold.
 \begin{itemize}
 \item[(i)]A codimension one distribution $\cD   \subset TS$  on $S$ is called {\it contact} if  the $2$-form
 $\o^\theta =  d \theta|_{\cD}$, where  $\theta$ is a defining form for $\cD$,   is non-degenerate.  The  defining form  $\theta$  is  called  {\it contact form}.
\item[(ii)]
The  {\it  Reeb vector field $Z = Z^\theta$}  of a  contact form $\theta$ is the unique vector field satisfying the conditions
\beq Z\lrcorner\, \theta=1,\qquad  Z\lrcorner\, d\theta=0\ .\eeq
\end{itemize}
\end{defi}
A Reeb vector field $Z = Z^\theta$ determines a direct sum decomposition $TS = \cD + \langle Z \rangle$ and preserves  $\theta$ (indeed,
 $\cL_Z \theta= d(Z\lrcorner \theta)+ Z\lrcorner d\theta=0$)  and the contact distribution  $ \cD = \ker \theta$.
The conformal class $[\theta]$  of the contact  forms  of $\cD$  is   globally  defined.
\begin{defi}
 A sub-Riemannian structure    $(\cD, h^\cD)$ (resp. a subconformal structure $(\cD, [h^\cD])$)  is called {\it of contact type} if the underlying distribution $\cD$ is contact.
\end{defi}

 \subsubsection{Almost,  partially integrable and integrable CR structures}
 \begin{defi}
 Let $\cD  \subset TS$ be a contact structure on a manifold $S$  and   $J $  a complex structure on $\cD$ (that is, a field $J \in \G(\End(\cD))$ of endomorphisms of $\cD$ with  $J^2 =  - \Id_{\cD}$). 
 The pair $(\cD, J)$  is called  {\it almost CR structure}. Moreover:
 \begin{itemize}[leftmargin = 20pt]
 \item[(i)]
An almost CR structure   $(\cD = \ker \theta, J)$  is called {\it partially integrable} if  the associate $2$-form $h^\theta = \o(\cdot, J \cdot)$, $\o = d \theta|_\cD$,
is symmetric.  It is called Levi  form and its conformal class $[h^\theta]$ is globally defined.
\item[(ii)]  A partially integrable   structure $(\cD, J)$ is called  {\it strongly pseudoconvex}  if    the   Levi form $h^{\theta}$  is  positive (or negative) definite.
\item[(iii)]  An {\it integrable CR structure}  $(\cD,J)$    is a  partially integrable CR structure   with identically   vanishing   Nijenhuis tensor $N_J \in \cD^* \otimes \cD^* \otimes \cD$,  defined by
$$N_J(X,Y)=[X,Y]-[JX,JY]+J([JX,Y]+[X,JY]) \ , \ \ X, Y \in \cD .$$
\end{itemize}
\end{defi}
\par
\medskip
\subsubsection{Correspondence between subconformal  and CR structures}
\begin{theo}[\cite{AGS}] \label{bijection}  Let $(S, \cD = \ker \theta)$ be a  contact manifold  with a globally defined contact form $\theta$. There exists a one-to-one correspondence between the following two sets.
\begin{itemize}[leftmargin = 15pt]
\item The pairs  $((\cD, J), B)$,   in which  $(\cD, J)$ is a strongly pseudoconvex almost  CR structure with positive Levi form  $h^ \theta >0$,  and  $B \in \End(\cD)$ is a $h^\theta$-symmetric and positive definite field of endomorphisms;
\item  The sub-Riemannian structures of contact type $(\cD, h)$.
\end{itemize}
Such a correspondence is  given by
\beq\label{mapsto}  ((\cD, J), B) \qquad \xleftrightarrow{\hspace*{0.7cm}}  \qquad (\cD, h \= h^{\theta} \circ B )\ .\eeq
\end{theo}
\begin{pf} It suffices to show that the correspondence  $((\cD, J), B)  \mapsto   (\cD, h \= h^{\theta} \circ B)$  has an inverse.  To prove this,  let $h$ be a sub-Riemannian metric   on $\cD$. Denote by   $K$ the     field  of endomorphisms  of $\cD$ defined by  $K \= h^{-1} \circ \o$, $\o \= d \theta|_\cD$. Note  that  $K$  is  $h$-skew-symmetric and, consequently,   that  $- K^2 = - K \circ K> 0$ is $h$-symmetric  and positive.   Consider the field of endomorphisms
$B \= (-K^2)^{-\frac{1}{2}}  > 0$, i.e.  the  inverse of the (unique) positive square root of $-K^2$. Then the field of endomorphisms   $J \= B K$ is a   complex structure on $\cD$ as the following calculation shows
\beq  J^2 =  B K B K =(-K^2)^{-\frac{1}{2}} K  (-K^2)^{-\frac{1}{2}} K
=  \left((-K^2)^{-\frac{1}{2}} \right)^{2} K^2 = - \Id_\cD\ .\eeq
 Moreover, since $ \o = h \circ K$  and $B$ is $h$-symmetric,  $h^\theta \= \o( \cdot ,  J \cdot) $ is symmetric:
\begin{multline}  h^\theta(X, Y) = \o( X,  J Y) = \o ( J^{-1} Y, X)  =   h\left(B^{-1} Y, X\right) =  h\left( Y, B^{-1} X\right)\\
=   - \o\left( Y, K^{-1} B^{-1} X\right) =  \o\left( Y,  J X\right) = h^\theta(Y, X) \ .\end{multline}
This means   that  $(\cD, J)$ is a partially integrable almost CR structure. Since $B$ is positive and $h^\theta = h \circ B^{-1}$, the almost CR structure  $(\cD, J)$ is strongly pseudoconvex. Hence, the correspondence $(\cD, h) \mapsto ((\cD, J), B)$ is the desired inverse map.
\end{pf}
By  the above   proof,   the  almost CR structure $(\cD, J)$
 depends only on the conformal class $[h]$ of the sub-Riemannian metric.   Hence Theorem \ref{bijection}  implies  the following  corollary.
\begin{cor}  \label{cor35} On a  contact manifold $(S, \cD)$ there is  a canonical one-to-one correspondence between the subconformal structures  $(\cD, [h])$ and the  pairs $((\cD, J), [B])$ formed by  a strongly pseudoconvex
almost  CR structure $(\cD, J)$ and a conformal class  $[B]$ of  fields of $\cD$-endomorphisms that are  positive definite with respect to the positive Levi forms of $(\cD, J)$.
\end{cor}
\par \medskip
\subsection{Sasaki and K\"ahler manifolds associated with    Robinson-Trautman bundles}
\subsubsection{Regular Sasaki manifolds}
Let   $(\cD = \ker \theta,J)$  be a strongly pseudoconvex  {\it integrable} CR structure on a contact manifold $(S, \cD = \ker \theta)$  equipped with  a fixed contact form $\theta$.
We recall that the Reeb vector field $Z = Z^{\theta}$ is transversal  to  $\cD$ and preserves $\theta$ and  $\cD$.
\begin{definition} \label{defsasaki} The   CR manifold  $(S, \cD = \ker \theta, J)$  is called {\it Sasaki}  if  the Reeb vector  field $Z = Z^\theta$ preserves   $J$, i.e.,   $\cL_Z J =0$. Such
  Sasaki manifold  is called {\it regular}  if   $Z$  generates a one-parameter  group $A= \exp(\bR Z)$ of diffeomorphisms
acting freely and properly on $S$.
\end{definition}
The  {\it  Sasaki metric}  of a Sasaki manifold  $(S, \cD = \ker \theta, J)$  is the  Riemannian  metric
\beq \label{Sasakimetric}  g^\theta \= \theta^2 +   \frac{1}{2}  h^\theta  \eeq
where  $h^\theta$ is  considered as a  degenerate  metric  on  $S$  with  kernel  $\langle Z\rangle$.
Note  that  the Reeb vector  field  $Z = Z^\theta$   preserves the   Sasaki metric $g^\theta$.\par
\smallskip
In the literature a Sasaki manifold  is often defined in  a different   way,  namely as    a Riemannian manifold equipped with  a  unit  Killing  vector  field $Z$
 satisfying appropriate  conditions.
The  following  proposition shows that  the two definitions  are   equivalent.\par
\begin{prop} \cite{ACHK} Let $( S, \cD = \ker \theta, J)$ be  a Sasaki manifold  with Reeb vector field $Z = Z^\theta$ and  Sasaki metric $g = g^\theta$.
Then  $Z$ is  a unit Killing vector field for $g$  and the pair  $(g, Z)$ satisfies the relations
  \be \label{Sasakiconditions}   \theta = g\circ Z  \qquad \text{and}\qquad   J = g^{-1} \circ d \theta |_{\cD} \ .\ee
  Conversely  any Riemannian manifold   $(S,g, Z)$  with a unit Killing vector field  $Z$ such that
  \begin{itemize}
  \item[(1)] $\theta = g\circ Z$ is a contact form;
  \item[(2)]   the  pair    $(\cD = \ker \theta,J)$,  with $J \= (g^{-1}\circ d\theta)_\cD =  \nabla^g Z|_{\cD} $,   is an integrable  pseudoconvex CR structure,
  \end{itemize}
   determines the Sasaki manifold  $(S, \cD = \ker \theta, J)$.
  \end{prop}
\medskip

\subsubsection{Correspondence between  Sasaki manifolds and quantisable K\"ahler manifolds}
\begin{definition} A K\"ahler manifold $(N, J^N, g^N)$  is called {\it quantisable} if there exists a principal $A$-bundle $\pi: S \to N$, with $A = S^1$ or $\mathbb{R}$,  with
a connection $1$-form $\theta: TS \to \bR$, whose curvature $d \theta$ is equal to  the K\"ahler form $\o^N = g^N(\cdot, J^N\cdot)$, more precisely     $d\theta= \pi^* \o^N$.
\end{definition}
 For a fixed  $A = S^1$ or $\bR$, the K\"ahler manifold $(N, J^N, g^N)$ is called {\it $A$-quantisable} if it satisfies the above condition assuming  the group is  $A$. 
By \cite{ACHK}*{Prop. 1.2}, $N$ is $S^1$-quantisable if and only if the (\v Cech) cohomology class  $[\o^N] \in H^2(M, \bR)$ is integral. It is $\bR$-quantisable if and only if $[\o^N] = 0$.  
\par
\smallskip
In the next theorem  we establish  a  natural correspondence  between the  regular Sasaki  manifolds and the quantisable K\"ahler manifolds.
\begin{theo} \cite{ACHK} \label{Theorem45}
Let $(S,\cD = \ker \theta,J)$ be a regular Sasaki manifold with  Reeb vector  field  $Z = Z^\theta$ and denote by $A = \{e^{tZ}\} \simeq \bR$ or $S^1$,   the  group of diffeomorphisms generated by $Z$.\par
Then    $\pi:  S \to N= S/A$  is a principal $A$-bundle and $\theta$ is  a  connection $1$-form for such a bundle.  Moreover, the  $A$-invariant  complex structure $J$ and the  $2$-form  $\omega = d\theta|_{\cD}$  on $\cD$  project   onto  an  integrable  complex structure $J^N$ and  a  symplectic form  $\o^N$  on $N$, respectively,   such that  $(N, J^N,  g^N = \o^N \circ J)$ is a  quantisable K\"ahler  manifold.\par
Conversely, if $(N,J^N, g^N)$ is  a quantisable K\"ahler manifold  and  $\pi: S \to N$ is a principle $A$-bundle with $A = \bR$,  $S^1$ and  a connection $1$-form $\theta$ such that
$d\theta= \pi^* \o^N$,   then  $(S, \cD = \ker \theta,  J)$ is a regular Sasaki manifold.
\end{theo}
\begin{pf} We only need to prove the second claim. For this, we observe that the equality  $d\theta= \pi^* \o^N$ implies  that  $\theta$ is a contact form,  whose associated
 Reeb vector field $Z = Z^\theta$ coincides with   the fundamental vector field   of  the principal bundle.  Let $J$ be  the   field of endomorphisms of $\cD$ defined by
   $J_u = (\pi_*|_{\cD_u})^{-1}(J^N_{\pi(u)})$, $u \in S$. Then $(\cD, J)$ is a $Z$-invariant strongly pseudoconvex integrable CR structure and  $(S, \cD = \ker \theta, J)$ is a regular Sasaki manifold.
\end{pf}
Given a regular Sasaki manifold  $(S, \cD = \ker \theta, J)$  the  {\it associated K\"ahler manifold} is  the quantisable K\"ahler manifold  $(N = S/A, J^N, g^N)$ defined in the above theorem.
\par
\bigskip
\subsubsection{Twisting RT bundles associated with    Sasaki and K\"ahler manifolds} \label{section323}
The following lemma establishes a fundamental  relation between      twisting RT structures and  the  subconformal structures of contact type.
\begin{lem} Let $(\pi: M \to S, (\cW, [h]))$ be an RT bundle of even dimension.
The  RT  structure $(\cW, [h])$ on the bundle  $\pi: M\to S$  is twisting  if and only if  the corresponding subconformal structure $(\cD, [h^\cD])$ on the base manifold  $S$ is of contact type.
\end{lem}
\begin{pf}  Let us denote by $\theta$  a  defining $1$-form  for the distribution $\cD = \pi_*(\cW) \subset TS$ of the subconformal structure on    $S$ and   by   $\q \= \pi^* \theta$  the corresponding defining $1$-form for   $\cW$ on $M$.  
The claim follows immediately from the fact  that $d\theta|_\cD$ is  non-degenerate  if and only if $\dim \ker d\q|_{\cW_x} = 1$ for any  $x \in M$, i.e.  $(\cW, [h])$ is twisting.
\end{pf}
This lemma and the previous discussion about subconformal structures of contact type motivates the following
\begin{defi} Let $(\pi: M \to S, (\cW, [h]))$ be a twisting RT bundle  of even dimension $n = 2 k +2$ and  $(\cD, [h^\cD])$ the corresponding  subconformal structure of contact type on $S$.
 Let also $((\cD, J), [B])$   be  the pair given by a strongly pseudoconvex almost CR structure and a conformal class $[B]$ of   positive definite endomorphisms, which  corresponds to   $(\cD, [h^\cD])$  by Corollary  \ref{cor35}.  The RT bundle    $(\pi: M \to S,(\cW, [h]))$   is called  {\it of K\"ahler-Sasaki type} if:
\begin{itemize}[leftmargin = 25pt]
\item[(a)]   there exists a global   contact form $\theta$ for the contact distribution $\cD = \ker \theta \subset TS$;
\item[(b)]  $[B] = [\Id_\cD]$, i.e.   $[h^\cD] = [h^\theta]$ is the conformal class of the positive Levi forms of $\cD$;
\item[(c)]   $(S, \cD = \ker \theta, J)$  is a  regular  Sasaki manifold, i.e. it  is principle $A$-bundle over  a $2k$-dimensional quantisable K\"ahler manifold $(N, J^N, g^N)$.
\end{itemize}
The  compatible  Lorentzian  metrics   of  $(M, (\cW, [h]))$   are   called {\it of  K\"ahler-Sasaki type}.
\end{defi}
\par
\bigskip
 \subsection{Generalised electromagnetic plane waves}
Let $(M, g)$ be an orientable Lorentz\-ian manifold with  canonical volume form $\operatorname{vol}_g$. For any $0 \leq p \leq n$,  the usual Hodge-${\ast}$ operator  ${\ast}: \O^p(M) \to \O^{n-p}(M)$ is defined by 
 \beq \a \wedge {\ast} \b = g(\a, \b) \operatorname{vol}_g\ ,\hskip 2 cm  \a \in \O^p(M)\ ,\ \ \b \in \O^{n-p}(M)\ . \eeq
 As we mentioned in the introduction, if $\dim M = 4$,   an  electromagnetic  plane wave on $(M, g)$ is  a  decomposable   $2$-form  $F = \q \wedge \operatorname{e}^*$, determined by a null $1$-form $\q$ and  a  $g$-orthogonal   space-like  $1$-form   $\operatorname{e}^*$, which is harmonic, that is such that $d F = 0 = d(\ast F)$. 
 \par
 \smallskip
According  to  Trautman (\cite{Tr2}),  there  exist  two natural possible ways to generalise the definition of  electromagnetic   plane wave for the Lorentzian manifolds  of higher dimension $n > 4$. The first way is to consider only  manifolds of even dimension $n = 2k$ and  assume that  a ``generalised  plane  wave''   is any  harmonic $k$-form,   which  is locally  of the form  $F = \q \wedge e^1 \wedge \cdots e^{k-1}$ for a $1$-form $\q$  with  $\ps \= g^{-1} \circ \q$  null and  $\{e^1, \cdots , e^{k-1}\}$  a set of    $g$-orthonormal  linearly independent space-like  1-forms. A second alternatively way (which can be used  for Lorentzian manifold of arbitrary dimensions) is just to use the    definition  used for the  $4$-manifolds. However,  it seems  that in higher dimensions such second definition  does no longer implies that the one-dimensional  distribution $\cK= \ker F \cap \ker ({\ast} F)$ is  generated by a  \sf vector field $\ps$ (\cite{Or}; see also \cite{SOMA}
 ). We think that  either  this or some other  property should be included as a part of the definition. It is   an issue that we leave to a future work.  Here we follow just the   first  way of generalising the notion of  electromagnetic plane wave and  we adopt the following
\begin{defi} A {\it generalised electromagnetic plane wave}  on a  Lorentzian manifold $(M, g)$ of dimension $\dim M = 2k$  is a harmonic   $k$-form $F$, which is a wedge product $F =  \q \wedge  \a$
of a null $1$-form $\q$ and a  $(k-1)$-form  $\a$  with the property  that any null vector of the distribution  $\cW_F = \ker \q$ is also  in  $ \ker \a$.
\end{defi}  
\begin{prop} \label{waves} Let $F = \q \wedge \a$ be a generalised  electromagnetic plane wave  on an oriented  Lorentzian $(n = 2k)$-manifold $(M, g)$ with $n \geq 4$ and $\cW_F = \ker \q$. 
Then:
\begin{itemize} 
\item[(1)] Also the  dual $k$-form ${\ast} F$ is  a generalised electromagnetic plane wave, hence of the form ${\ast} F = \q \wedge \b$;  
\item[(2)] Let $\cK_F = \langle \ps_o \rangle$ where $\ps_o = g^{-1} \circ \q$. Then  $ \cK_F  =\ker F \cap \ker ({\ast} F)$ and
\beq \label{lieder} \cL_{\ps_o} F = 0\qquad \text{and}\qquad \cL_{\ps_o}({\ast} F) = 0\ .\eeq
\item[(3)]  For any  vector field $\ps   \in \cK_F  = \langle \ps_o \rangle$ 
\begin{align}
\label{1} & \cL_\ps\q  = -  f \q \qquad \text{and thus}\qquad \cL_\ps \cW_F \subset \cW_F\ ,\\
\label{2}  &  \cL_{\ps}\a  =  f \a + \q \wedge \g \ ,\qquad  \cL_{\ps}\b  =   f \b + \q \wedge \g'\ .
 \end{align}
where   $f$ is a function and $\g, \g' $  are  $(k-2)$-forms with $\cK_F \subset \ker \g \cap \ker \g'$. 
\end{itemize}
\end{prop}
\begin{pf} (1) Since $\q$ is null, there exist a null vector field $\ps_o$ such that $\q = \ps_o^\flat = g(\ps_o, \cdot)$ and a null vector field $\qs_o$ such that $\q(\qs_o) =  g(\ps_o, \qs_o) = 1$.  
Consider a (local) frame field $(\ps_o, e_1, \ldots, e_{2k-2}, \qs_o)$
where    $( e_1, \ldots, e_{2k-2})$  is a $g$-orthonormal basis for $\cE = \langle \ps_o, \qs_o \rangle^\perp$ and let $(\q', e^1, \ldots, e^{k-2}, \q)$ be its  dual coframe field. Note that $\cW_F = \ker \q = \langle \ps_o, e_1, \ldots, e_{2k-2}\rangle$ and $\ker g_{\cW_F} = \langle \ps_o \rangle$.  Moreover, since  
 $ \langle \ps_o \rangle = \ker g_{\cW_F}  \subset \ker \a$ and  $\a$ is determined up to terms of the form  $\q \wedge \g$,  we may assume that 
$$\a = \sum_{i_1 < \ldots < i_{k-1}} \a_{i_1 \ldots i_{k-1}} e^{i_1} \wedge \ldots e^{i_{k-1}}\ .$$ 
Hence  ${\ast} F = \q \wedge \b$ with  $\b =  {\ast}_\cE \a$,  where we denote by ${\ast}_\cE$ the  $\ast$-Hodge operator of the Euclidean space $\cE$.  In other words, ${\ast} F = \q \wedge \b$ for a  $(k-1)$-form $\b$ such that   $ (\a \wedge \b)_{\cE}  = \operatorname{vol}_\cE$ and it is therefore a generalised plane wave. \par
(2)  Since $d F = 0 = d({\ast} F)$ and $\ps_o \in \ker F \cap \ker ({\ast} F)$,  the Lie derivatives along $\ps_o$ of $F$ and ${\ast}F$ are trivial.  Furthermore, $\ps \in \ker F \cap \ker ({\ast} F)$ if and only if 
$\ps \lrcorner \q = 0$ and $\ps\lrcorner (\a \wedge \b) = 0$.  Since  $(\a \wedge \b)_{\cE}  = \operatorname{vol}_\cE$, this occurs if and only if $\ps\in \langle \ps_o \rangle  = \cK_F$.  \par
\smallskip
(3)  It is sufficent to prove the claim for  $\ps = \ps_o$. From  (2) we have that  
  \beq \label{deriv} 
   0 = \cL_{\ps_o}F = \cL_{\ps_o}   \q \wedge\a  + \q \wedge \cL_{\ps_o} \a \ , \qquad 0 =   \cL_{\ps_o}({\ast} F) = \cL_{\ps_o}   \q \wedge\b  + \q \wedge \cL_{\ps_o} \b\
  \eeq
We  expand $d \q$ as  
  \beq
   d\q =  \q \wedge e^* + f \q \wedge \q'\ + \sum_{i \leq j}  \nu_{ij} e^i \wedge \e^j +  e'{}^* \wedge \q'  
  \eeq
  for some $e^*, e'{}^*  \in \cE^* = \langle e^1, \ldots, e^{2k-2}\rangle$. Hence
  $$\cL_{\ps_o}F = - f   \q \wedge\a  - e'{}^* \wedge \a + \q \wedge \wt \a  \ , \qquad  \cL_{\ps_o}({\ast} F) = - f \q \wedge\b  -e'{}^* \wedge \b + \q \wedge  \wt \b $$
 where  $\wt \a, \wt \b$ are  the  $(k-1)$-forms  such that 
 $$\cL_{\ps_o} \a = \wt \a \mod \left\langle \q \wedge \g\ ,\ \g \in \L^{k-2}  \cE^* \right\rangle \ ,\qquad \cL_{\ps_o} \b = \wt \b \mod \left\langle  \q \wedge \g, \  \g \in \L^{k-2} \cE^*\right\rangle\ .$$
Therefore the vanishing    $ \cL_{\ps_o}F  = \cL_{\ps_o}({\ast} F) = 0$ implies that 
  $$  \cL_\ps \a  =   f \a\mod \left\langle \q \wedge \g\ ,\ \g \in \L^{k-2} \cE^* \right\rangle \ ,\quad  \cL_\ps \b  =   f \b\mod \left\langle \q \wedge \g\ ,\ \g \in \L^{k-2} \cE^* \right\rangle$$
  and  $e'{}^* \wedge \a  =  e'{}^* \wedge \b= 0$.  Since $(\a \wedge \b)_{\cE} =\operatorname{vol}_\cE$, this  implies   $e'^* = 0$ and   
  $ \cL_{\ps_o} \q  = -  f \q$.
\end{pf}
\begin{defi}  Let  $(M, g)$ be a Lorentzian $2k$-manifold. A {\it flag structure} on $M$ is  a pair $(\cK, \cW = \ker \q)$, determined by a  null $1$-form  $\q$ and a   one-dimensional distribution  $\cK = \ker g_\cW$. Any   generalised electromagnetic plane wave    $F = \q \wedge\a$ on  $(M, g)$  determines the  flag structure 
$$(\cK_F\= \ker F \cap \ker( {\ast} F) = \ker g_{\cW_F}\ ,\ \cW_F \= \ker \q)$$
 which we call the   {\it  flag   structure} of $F$.  
\end{defi}
  Proposition \ref{waves}   implies the following
     \begin{cor} \label{thecor} Let $(M, g)$ be a  manifold of $\dim M = 2k$, equipped with  a flag structue  $(\cK = \ker g_\cW, \cW = \ker \q)$. A necessary condition for the existence of a 
     generalised  plane wave $F$  having $(\cK, \cW)$ as its  flag structure  is the existence of  two $(k-1)$-forms $\a, \b$  satisfying the following  conditions: 
     \begin{itemize}
     \item[(a)]  $\cK\subset \ker \a_\cW \cap \ker \b_\cW$; 
     \item[(b)] $\q \wedge \b = \ast (\q \wedge \a)$ ;  
     \item[(c)] any vector field  $\ps \in \cK$  preserves $\cW$  and there is a function $f$ such  that 
$$   \cL_\ps \a =  f \a + \q \wedge \g \ ,\qquad  \cL_{\ps}\b  =  f \b + \q \wedge \g'$$
for some $\g, \g' $ such that $\ker g_{\cW_F} \subset \ker \g \cap \ker \g'$. 
\end{itemize}
    \end{cor}
  \begin{rem} \label{rem211}  If the flag structure     $(\cK, \cW = \ker \q)$ is  determined by a \sf structure, then the  necessary conditions of  Corollary  \ref{thecor}
  are  satisfied.  If in addition $\dim M = 4$, such necessary  conditions  are equivalent to say that $(\cK, \cW)$ is the pair  determined by a \sf structure.
  \end{rem}
     \par
     \medskip
\subsection{ The Robinson Theorem for  Lorentzian manifolds of K\"ahler-Sasaki type}
 \begin{theo}\label{Robinson}
 Let $(M, g)$ be a $(n=2k)$-dimensional  Lorentzian manifold of K\"ahler-Sasaki type with associated \sf structure $(\cW, [h = g_\cW])$.  Then locally  there exists a  non-trivial generalised  electromagnetic plane wave $F$ with  flag structure   $(\cK_h = \ker h, \cW)$.
 \end{theo}
 \begin{pf}   Let  $\pi: M \to S$ and  $\pi^S: S \to N$ be the  principal bundles over  the regular Sasaki manifold $(S = M/A, \cD = \ker \theta, J)$ and  the quantisable K\"ahler manifold $(N = S/A', J, g^N)$.
  Let     also    $(z^1, \cdots z^{k-1})$  be local holomorphic  coordinates on $N$ and  denote by  $\f = \f(z,\bar z)$  and  $\omega = i\p \bar \p \varphi$ a potential and the K\"ahler form of $N$.
 The connection  $1$-form  $\theta = -  \frac{1}{2} du  -  \frac{1}{2} d^c \varphi=  -  \frac{1}{2} du -\frac{i}{2}(\overline{\partial} \varphi - \partial \varphi)$  has curvature $ d \theta = i \partial  \overline{\partial}\varphi  =  \omega$. 
 \par
Since all arguments are local, we may assume that  $M$  is a trivial bundle   $\pi: M = A \times S\to S$  with $A = \bR$ or $S^1$. We denote by  $t$ its fiber coordinate.  Any  admissible metric  on $M$  has locally the  form
 $$   g =  g^N + \q \vee \eta  $$
 where $g^N$ and  $\q$ are   the pull-backs  of the metric of $N$ and of  the  1-form $\theta$ of $S$,  respectively,  and  $\eta \in \Omega^1(M)$ is a 1-form which can be locally written as 
 $  \eta =  \alpha  dt + \gamma^i \varphi_{i\bar{j}} dz^j  + \overline{ \gamma^i \varphi_{i\bar{j}} dz^j } + \beta  \vartheta$ for some real functions  $\alpha, \beta$  and complex functions  $\gamma^i$.\par
 \smallskip
 Let $\cF$ and $F$  be   the  complex $k$-form  $\cF = \vartheta \wedge  dz^1 \wedge \cdots \wedge dz^{k-1} $ and 
  the real $k$-form
 $$F = \mathrm{Re}\cF = \vartheta \wedge (  dz^1 \wedge \cdots \wedge dz^{k-1} +\overline{dz^1 \wedge \cdots \wedge dz^{k-1}}         )\ ,$$
respectively.  By Lemma \ref{concludinglemma} below, $F$ is harmonic and thus it is  a  generalised electromagnetic  plane wave. \end{pf}
 \begin{lemma}    \label{concludinglemma} The  complex  form  $\cF$ is closed and  coclosed.
 \end{lemma}
  \begin{pf} For the closedness, just observe that    $d \cF = d (\vartheta \wedge dz^1 \wedge \cdots \wedge dz^{k-1} )
   = \omega \wedge dz^1 \cdots \wedge  dz^{k-1} =0$.
The co-closedness follows from the fact that ${\ast} \cF =  \pm  i^{k-1} \cF$ (the sign  depending on the orientation) since at  each  point  $x \in M$ the space $\langle d\theta, dz^1, \ldots, dz^{k-1}\rangle|_x$ 
is an isotropic subspace of $(T^\bC_x M, g_x)$.  \end{pf}
The  proof has  the following  consequence. 
\begin{cor}   Assume that  the Lorentzian manifold $(M, g)$   of K\"ahler-Sasaki type  is globally trivial $M = \bR \times S$ and that  the corresponding  K\"ahler manifold  $(N, J^N, g^N)$ has a    {\rm global}  holomorphic  volume  form. Then there exists on $M$  a nowhere vanishing {\rm globally defined}  generalised plane wave.  \end{cor}
\par
\medskip 
\section{Einstein metrics on  Robinson-Trautman bundles}\label{section4}
 \label{notation}
Throughout this  section $\pi: M = S \times \bR \to S$ is  a trivial   $n$-dimensional principal bundle  with structure group $A = \bR$ over a regular Sasaki manifold  $(S, \cD = \ker \theta, J)$ fibering over a quantisable  K\"ahler manifold $(N = S/A', J,  g^N)$
  with  structure group $A' = \bR$ or $S^1$.  The  assumption made here   that  the $\bR$-bundle $\pi: M \to S$ is  trivial  is mostly chosen for the sake of  simplicity. In fact, most parts of the
  following discussion remain valid under the weaker hypothesis that  such a  bundle  is  equipped with  a  flat connection.\par
  \medskip
  The following notation   is used.
  \begin{itemize}[leftmargin= 20pt]
  \item[--]  $\q := \pi^*(\theta)$  is  the pull back  of  the  contact form     of $S$ and
    $\cW = \ker \q$  is  the    corresponding   kernel   distribution  on $M$;
     \item[--]  $[h_o]$  is the conformal class of the degenerate metric $h_o = (\pi\circ \pi^S)^*g_o$  on $\cW$;
  \item[--]   $\cH_o = TS \subset TM$ is  the standard flat connection of the trivial $\bR$-bundle $M = S \times \bR$;
  \item[--]   $\ps_o = \frac{\p}{\p t}  \in \gX(M)$ and  $\qs^S_o \in \gX(S)$  are the   fundamental  vector fields  of the principal bundles  $\pi: M \to S$ and  $\pi^S: S \to N$, respectively;
  \item[--]   $\qs_o \in \gX(M)$  is the  horizontal lift of $\qs_o^S$ on $M$   with respect  to the  flat connection;
 \item[--]  for any  vector field $X$ on  $N$  we denote by:
\begin{itemize}[leftmargin = 20pt]
\item[$\cdot$] $X^{(S)} \= X^{\operatorname{h}}$  the  horizontal lift of $X$ in $\cD  = \ker \theta \subset TS$;
\item[$\cdot$] $\wh X = X^{(S)\operatorname{h}}$  the $\bR$-invariant horizontal  lift   of  $X^{(S)} $ in  $\cH_o \subset TM$.
\end{itemize}
  \end{itemize}
  Note that $(\pi: M  = S \times \bR \to S, (\cW, [h_o]))$ is an RT bundle of K\"ahler-Sasaki type, the kernel distribution   is   $\cK_h:= \ker h_{o} = \langle \ps_o \rangle$ and the pair
 $(\cW'_o = \cW \cap \cH_o, \qs_o)$ is a rigging for the \sf structure $(\cW, [h_o])$.
 Moreover,
  \begin{itemize}[leftmargin = 20pt]
 \item[(1)] $\cW'_o\subset TM$  is  the  horizontal lift (with respect to the flat connection) of $ \cD = \ker \theta \subset TS$;
  \item[(2)]   $\qs^S_o$ is  the Reeb vector field of  the contact  form $\theta$;
  \item[(3)] $d \theta|_{\cD}  = (\pi^S)^*(\o)$ where $\o \=  g_o( \cdot,  J\cdot)$ is the K\"ahler form of $(N, J, g_o)$;
  \item[(4)] for any pair of vector fields $X, Y$ on $N$, the corresponding  horizontal lifts $\wh X, \wh Y \in \gX(M)$ satisfy the  relations
\beq \label{62}  [\wh X, \wh Y] =   \wh{[X, Y]}  - g_o(X, JY) \qs_o\ ,\qquad [\wh X, \ps_o] = [\wh X, \qs_o] =  [\ps_o, \qs_o] = 0\ .
\eeq
  \end{itemize}
  We conclude this section by observing that if   $(E_i)$ is a (local) frame field for the K\"ahler manifold $(N, J, g_o)$, the corresponding  horizontal  lifts $\wh E_i$  on $M$  form  a frame   field for  $\cW'_o$  and  the tuple of vector fields  $(\ps_o, \wh E_i, \qs)$ is a (local) frame field on $M$.   We  denote by $(\ps_o^*,\wh E^i,  \qs_o^* = \q)$ the corresponding dual frame field.
\par
\medskip
In this section   we  determine a new family of    Lorentzian Einstein metrics on such RT bundle  belonging  to the following  special class of   compatible metrics.
 \smallskip
 \subsection{Firmly  compatible metrics} \label{severe}
By  Theorem \ref{compatibility}, any compatible metric  on $(M = S \times \bR, \cW, [h])$  is {\it locally} determined by a triple $(h, \ps, (\cW', \qs))$,  given by a degenerate metric $h = \s h_o \in [h_o]$, a vector field $\ps \in \cK_h$ and a rigging $(\cW', \qs)$.  With no loss of generality,   we    assume that $\ps = \ps_o = \frac{\p}{\p t}$ and 
the complementary subdistribution  is  $\cW' = \cW'_o$. \par
\smallskip
We will consider only   {\it globally defined} compatible metrics,  associated with  triples $(h, \ps$, $(\cW', \qs))$  where   $\qs$ is a {\it global} vector field of the form
\beq \label{52bis}   \qs \=a \qs_o + b \ps_o + E\ ,\qquad a \neq 0\ .\eeq
for some (global) functions $a$, $b$ and   vector field $E \in \cW_o'$.  The metrics that are associated with  a vector field $\qs$, for   which the coefficient  $a$ is  constant, are called {\it  strongly compatible}.  Our main results  deal with the following even more restricted  class of  compatible metrics.
\begin{definition} A compatible Lorentzian  metric $g$ on $M$  is called {\it   firmly  compatible} if it is determined by  a triple $(h, \ps, (\cW', \qs))$ as above,  in which the   coefficient   $a$ of $\qs$  is  constant  and the vector field $E$ is zero.
\end{definition}
Let $g$ be a  firmly compatible metric determined by a triple $(h = \s h_o, \ps_o, (\cW'_o, \qs))$. Let also  $(E_i)$ be  a local frame field of  the K\"ahler manifold
$(N, J, g_o)$ and  let $(\ps_o, \wh E_i, \qs)$, $(\ps_o^*,\wh E^i,  \qs_o^* = \q)$ be the corresponding  frame and dual coframe  fields  on $M$.   Then  by \eqref{buona}
any  firmly compatible metric on $M$  has the form
 \begin{multline} \label{buona-sev} g = \s h_{ij} \wh E^i\vee \wh E^j + \q \vee \left(\frac{2}{a}\ps^*_o  - \frac{2 b}{a^2} \q\right) =  \s\left( h_{ij} \wh E^i\vee \wh E^j +  \q \vee \left(\a\ps^*_o   + \b  \q\right) \right)\\
 \text{where}\ \a \= \frac{2}{a\s}\ ,\quad\quad \b \= - \frac{2 b}{a^2 \s}\ .\end{multline}
 Since  $a$ is  constant and any  homothetic rescaling of  the  metric $g_o$ on $N$ corresponds to an (inverse) homothetic rescaling of the Reeb vector field $\qs_o$  on $S$,
 with no loss of generality we may assume that  $a = 2$ and hence, by \eqref{buona-sev}, that
 \beq \a = \frac{1}{\s}\eeq
 and
 \beq \label{buona-sev1} g =     \s h_{ij} \wh E^i\vee \wh E^j +  \q \vee \left(\ps^*_o   + \wt \b  \q\right) =  \s  (\pi^S \circ \pi)^*(g_o)  +  \q \vee \left(\ps^*_o   + \wt \b  \q\right) \ \ \text{where} \ \  \wt \b \= \s \b = - \frac{b}{2}\ .\eeq
 \par
 \medskip
\subsection{Einstein  metrics of Taub-NUT type} In the next theorem we  describe   the  family of Einstein metrics we advertised in the Introduction.
 \begin{theo} \label{main}  Let $\pi: M =  S \times \bR \to S$ and  $\pi^S: S \to N$ be as above and assume that the quantisable K\"ahler manifold $(N, J, g_o)$  is Einstein  
 (\footnote{We recall  that   this occurs if and only if     $S$ is  Sasaki-Einstein (see e.g. \cite{BG}*{Ch.11}).}) with Einstein constant $\L_o$. Furthermore, for any  triple of real numbers $(\L, B, C)$  with   $C>0$,  let   $\s: \bR \to \bR$ and  $\wt \b: (0, + \infty) \to \bR$  be the functions defined by
 \beq \label{eqtheorem0}  \s(t) \=   \frac{1}{16 C}  t^2 + C\ ,\hskip 10 cm \eeq
\beq \label{eqtheorem}
 \wt \b(t) \=   \frac{t}{(t^2+16C^2)^{\frac{n}{2}-1}} \left(B - \int^t_1 \frac{\left(16 C^2+s^2\right)^{\frac{n}{2}-1}  \left(16 C \Lambda _0- \Lambda \left(16 C^2+s^2\right)\right)}{4 s^2}d s\right)\eeq
  Then $\wt \b$ is a rational function  admitting a unique smooth extension over $\bR$ and the corresponding    firmly compatible  metric $g$ on $(M = S \times \bR, (\cW, [h_o]))$
  \beq \label{metrictheorem} g  =   \s  (\pi^S \circ \pi)^*(g_o)  + \q_o \vee ( \ps_o^* +   \wt \b  \q_o )\ .\eeq
is  Einstein   with Einstein constant $\L$.\par
Conversely, if  $g$ is a metric on $M = S \times \bR$, which is 
\begin{itemize}
\item[(a)] firmly compatible, hence of the form  \eqref{metrictheorem} for some functions $\s, \wt \b = \s \b$ and
\item[(b)]  with  function   $\s, \wt \b$   depending just on the  coordinate $t \in \bR$, 
\end{itemize}
 then $g$  is Einstein with Einstein constant $\L$ if and only if  $\s$ and  $\wt \b$ are the functions defined in \eqref{eqtheorem0} and \eqref{eqtheorem} for some choice of the constants  $B$ and $C>0$.
  \end{theo}
 \begin{pf}
Consider a (local) frame field $(\ps_o, \wh E_i, \qs)$ and the dual coframe field   $(\ps_o^*, \wh E^i, \qs_o^* = \q)$ on $M$ as described in \S \ref{severe}. We set
$$g_{ij} \= g_o(E_i, E_j)\ ,\ \  \o_{ij} \= g_o(E_i, J E_j)\ ,\ \  J_i^j = g^{jk} \o_{ki}\ ,\ \ [E_i, E_j] = c_{ij}^k E_k\ ,\ \  \wt \b \=  \b \s\ . $$
 According to   this notation, we  have   $J(E_j) = J_j^\ell E_\ell$.  As we mentioned above,
   in terms of the above coframe field,  any   firmly compatible  metric  with $a \equiv 2$  has the form \eqref{buona-sev1}. A  tedious (but straightforward) computation based on Koszul's formula (see \S Appendix \ref{Tedious} for details)  shows that  the Christoffel symbols  of the Levi-Civita connection in this frame field, i.e.  the functions $\GGa A B C$  defined by  
\begin{align*} & \n_{E_i} E_j = \GGa i j k E_k + \GGa i j {\ps_o} \ps_o + \GGa i j {\qs_o} \qs_o\ ,&&
 \n_{E_i} \ps_o = \GGa i {\ps_o} k E_k + \GGa i {\ps_o} {\ps_o} \ps_o + \GGa i {\ps_o} {\qs_o} \qs_o\ ,\\
& \n_{E_i} \qs_o = \GGa i {\qs_o} k E_k + \GGa i {\qs_o} {\ps_o} \ps_o + \GGa i {\qs_o} {\qs_o} \qs_o\ ,&&
 \n_{\ps_o} E_j = \GGa {\ps_o} j k E_k + \GGa {\ps_o} j {\ps_o} \ps_o + \GGa {\ps_o} j {\qs_o} \qs_o\ ,\\
 &  \n_{{\ps_o}} \ps_o = \GGa {\ps_o} {\ps_o} k E_k + \GGa {\ps_o} {\ps_o} {\ps_o} \ps_o + \GGa {\ps_o} {\ps_o} {\qs_o} \qs_o\ ,&&
 \n_{{\ps_o}} \qs_o = \GGa  {\ps_o} {\qs_o} k E_k + \GGa {\ps_o} {\qs_o} {\ps_o} \ps_o + \GGa {\ps_o} {\qs_o} {\qs_o} \qs_o\ ,\\
& \n_{\qs_o} E_j = \GGa {\qs_o} j k E_k + \GGa {\qs_o} j {\ps_o} \ps_o + \GGa {\qs_o} j {\qs_o} \qs_o\ ,&&
 \n_{{\qs_o}} \ps_o = \GGa {\qs_o} {\ps_o} k E_k + \GGa {\qs_o} {\ps_o} {\ps_o} \ps_o + \GGa {\qs_o} {\ps_o} {\qs_o} \qs_o\ ,\\
& \n_{{\qs_o}} \qs_o = \GGa  {\qs_o} {\qs_o} k E_k + \GGa {\qs_o} {\qs_o} {\ps_o} \ps_o + \GGa {\qs_o} {\qs_o} {\qs_o} \qs_o\ ,
\end{align*}
are   equal to  
 \begin{align}
\nonumber & \GGa i j m =
  g^{mk} g_o(\n^o_{E_i}   E_j, E_k)
+ \frac{1}{2 \s }  \wh E_i(\s) \d_j^m  + \frac{1}{2 \s }  \wh E_j(\s) \d_i^m
 -   \frac{1}{2 \s} g_{ij}     g^{mk} \wh E_k(\s)
\ ,\\
 \label{6.12} & \GGa i j {\ps_o}  =    g_{ij}  \left( - \qs_o(\s) +  2  \wt \b  \ps_o(\s) \right)\ ,\quad
 \GGa i j {\qs_o}  = - \frac{  \o_{ij}}{2} -  g_{ij}   \ps_o(\s)\ ,\\[12 pt]
 & \GGa i {\ps_o} m = \GGa {\ps_o} i m =  -   \frac{J^m_{i} }{4 \s}  +  \ps_o(\s) \frac{\d^m_i}{2 \s} \ ,\qquad  \GGa i {\ps_o} {\ps_o} =   \GGa  {\ps_o} i {\ps_o} =
 \GGa i {\ps_o} {\qs_o} =  \GGa {\ps_o} i  {\qs_o} = 0\ ,\\[12 pt]
&  \GGa i {\qs_o} m =  \GGa {\qs_o} i m = -  \wt \b  \frac{J_{i}^m }{2 \s}
+ \frac{ \qs_o(\s)\d_{i}^m}{2 \s} \ ,\ \
 \GGa i {\qs_o} {\ps_o} = \GGa  {\qs_o} i {\ps_o} =   \wh E_i(\wt \b) \ ,
\end{align}
\begin{align}
 &  \GGa {\ps_o}  {\ps_o}  m  =   \GGa {\ps_o}   {\ps_o} {\ps_o}  =   \GGa {\ps_o}   {\ps_o} {\qs_o}  =  0 \ ,\qquad   \GGa i {\qs_o} {\qs_o} =  \GGa  {\qs_o} i {\qs_o} = 0\ , \\[12 pt]
   &  \GGa {\ps_o}  {\qs_o}  m  = \GGa {\qs_o}  {\ps_o}  m  =  0 \ ,\qquad
 \GGa {\ps_o}   {\qs_o} {\ps_o}  =  \GGa {\qs_o}   {\ps_o} {\ps_o}  = \ps_o(\wt \b) \ ,\qquad \GGa {\ps_o}   {\qs_o} {\qs_o} =  \GGa {\qs_o}   {\ps_o} {\qs_o}   = 0  \ ,
\\[12pt]
 \label{6.17}  & \GGa {\qs_o}  {\qs_o}  m = -\frac{g^{mk}}{2\s}\wh E_k(\wt \b)\ ,\qquad  \GGa {\qs_o}   {\qs_o} {\ps_o}  =  \qs_o(\wt \b ) +  2 \wt \b  \ps_o(\wt \b ) \ ,\qquad
 \GGa {\qs_o}   {\qs_o} {\qs_o}  =  -  \ps_o( \wt \b ) \ .
  \end{align}
 Using these expressions, we may directly compute   the components $\RR A B C D$  of the   Riemann curvature tensor   (\footnote{Following   \cite{KN},  we  define   $R$  by the formula  $R_{X Y} Z :=  \n_{X} \n_{Y} Z -  \n_{Y} \n_{X} Z - \n_{[X, Y]} Z $.})
 in the   frame $(X_A) = (\ps_o,\wh E_1,\ldots, \wh E_{n-2},  \qs_o)$.  Note  that if   $X_A$,  $X_B$  are  commuting vector fields  of  the frame, then the corresponding components $\RR A B C D$ are given by the formula
  \beq\nonumber  \RR A B C D = X_A (\GGa B C D) -  X_B (\GGa A C D) - \GGa A C  F \GGa B F D +    \GGa B C F \GGa A F  D . \eeq
 By \eqref{62},  the only non-commuting pairs  in the  frame   are  those with    $X_A = \wh E_i$ and $X_B = \wh E_j$.  The  corresponding  components  $\RR i j C D$ are given by
  \beq\nonumber  \RR i j C D =  \wh E_i (\GGa j C D) -  \wh E_j (\GGa i C D) -  \GGa i C  F \GGa  j  F D +   \GGa j C F \GGa  i F D- c^k_{ij} \GGa k C D+ \o_{ij} \GGa {\qs_o} C D\ . \eeq
     Using these two  expressions, we can determine all  components  $\Ric_{AB} = \RR D A B D $  of the Ricci tensor and   write down  the Einstein equations $\Ric_{AB}  = \Lambda g_{AB}$. We list these equations below. In those expressions,  the terms that are struck out  are those which are immediately seen to be   $0$ on the base of  the above  expressions for the Christoffel symbols $\GGa A B C$. We also use the  shorthand  notation
  \beq \label{hatricci} \wh \Ric_{ij} \=     \wh E_m (\GGa i j m) -  \wh E_i (\GGa m j m) -  \GGa m j  \ell \GGa i  \ell  m   +   \GGa i j \ell  \GGa  m \ell  m - c^r_{mi} \GGa r j m\ . \eeq
 Under the ansatz \eqref{ansatz},    $\wh \Ric_{ij}$ is  equal to the   pull-back  on $M$ of the Ricci tensor  $R_{ij}$ of the base manifold $(N, g_o)$ (see the observations after \eqref{ansatz} below). 
Using this notation, the Einstein equations for  a metric of the form \eqref{buona-sev1} take the form
 \begin{align}
\nonumber  & \Ric_{ij}\  {=} \ \RR m i j m + \RR {\ps_o} i j {\ps_o} + \RR {\qs_o} i j {\qs_o}  = \\
\nonumber  & \hskip 1 cm {=} \   \wh E_m (\GGa i j m) -  \wh E_i (\GGa m j m) -  \GGa m j  \ell \GGa i  \ell  m   -  \GGa m j  {\ps_o} \GGa i {\ps_o}  m   -  \GGa m j  {\qs_o} \GGa i {\qs_o}  m +\\
\nonumber  &\hskip 1.5 cm  +   \GGa i j \ell  \GGa m \ell  m  +   \GGa i j {\ps_o}  \GGa m {\ps_o}  m +   \GGa i j {\qs_o}  \GGa m {\qs_o}  m - c^r_{mi} \GGa r j m  + \o_{mi} \GGa {\qs_o} j m +\\
\nonumber &\hskip 1.5 cm+  \ps_o (\GGa i j {\ps_o}) -  \xcancel{\wh E_i (\GGa {\ps_o} j {\ps_o})} -  \GGa {\ps_o} j  \ell \GGa i \ell  {\ps_o}  - \xcancel{\GGa {\ps_o} j  {\ps_o} \GGa i {\ps_o} {\ps_o}}  - \xcancel{ \GGa {\ps_o} j  {\qs_o} \GGa i {\qs_o}  {\ps_o} } + \\
\nonumber &\hskip 1.5 cm +  \xcancel{ \GGa i j \ell \GGa \ell {\ps_o} {\ps_o}} +  \xcancel{\GGa i j {\ps_o} \GGa {\ps_o} {\ps_o} {\ps_o} }+  \GGa i j {\qs_o} \GGa  {\ps_o} {\qs_o} {\ps_o} + \\
\nonumber & \hskip 1.5 cm+ \qs_o (\GGa i j {\qs_o}) - \xcancel{ \wh E_i (\GGa {\qs_o} j {\qs_o}) }-  \GGa {\qs_o} j  \ell \GGa i  \ell {\qs_o} - \xcancel{ \GGa {\qs_o} j  {\ps_o} \GGa i {\ps_o}  {\qs_o}} - \xcancel{ \GGa {\qs_o} j  {\qs_o} \GGa i {\qs_o}  {\qs_o}} +  \\
\nonumber & \hskip 1.5 cm + \xcancel{ \GGa i j \ell \GGa  {\qs_o} \ell {\qs_o} } +   \GGa i j {\ps_o} \GGa {\qs_o}  {\ps_o} {\qs_o}  +   \GGa i j {\qs_o} \GGa {\qs_o} {\qs_o} {\qs_o} =\\[10pt]
\nonumber & \hskip 1 cm{=}\ \wh \Ric_{ij}  -  \GGa m j  {\ps_o} \GGa i  {\ps_o}  m   -  \GGa m j  {\qs_o} \GGa i  {\qs_o}  m + \GGa i j {\ps_o}  \GGa m {\ps_o}  m+   \GGa i j {\qs_o}  \GGa m  {\qs_o}  m    + \o_{mi} \GGa {\qs_o} j m +\\
\nonumber &\hskip 1.5 cm+  \ps_o (\GGa i j {\ps_o})  -  \GGa {\ps_o} j  \ell \GGa i \ell  {\ps_o}  +   \GGa i j {\qs_o} \GGa{\ps_o}  {\qs_o}  {\ps_o} + \qs_o (\GGa i j {\qs_o}) - \GGa {\qs_o} j  \ell \GGa i\ell  {\qs_o} +  \\
  \label{6.20} & \hskip 1.5 cm  +  \xcancel{  \GGa i j {\ps_o} \GGa {\qs_o} {\ps_o}  {\qs_o} }  +   \GGa i j {\qs_o} \GGa {\qs_o} {\qs_o} {\qs_o} =    \s \Lambda g_{ij}
 \end{align}
  \begin{align}
\nonumber & \Ric_{i\ps_o}   {=} \ \RR m i {\ps_o} m + \RR {\ps_o} i {\ps_o} {\ps_o} + \RR {\qs_o} i {\ps_o} {\qs_o}  = \\
\nonumber  & \hskip 1 cm {=} \   \wh E_m (\GGa i {\ps_o} m) -  \wh E_i (\GGa m {\ps_o} m) -  \GGa m {\ps_o}  \ell \GGa i  \ell  m   -  \xcancel{\GGa m {\ps_o}  {\ps_o} \GGa i  {\ps_o}  m }  - \xcancel{ \GGa m {\ps_o}  {\qs_o} \GGa i {\qs_o}  m }+\\
\nonumber  &\hskip 1.5 cm  +   \GGa i {\ps_o} \ell  \GGa m  \ell  m  +  \xcancel{ \GGa i {\ps_o} {\ps_o}  \GGa m  {\ps_o}  m } +   \xcancel{\GGa i {\ps_o} {\qs_o}  \GGa m  {\qs_o}  m }- c^r_{mi} \GGa r {\ps_o} m   + \xcancel{ \o_{mi} \GGa {\qs_o} {\ps_o} m }+\\
\nonumber &\hskip 1.5 cm+ \xcancel{ \ps_o (\GGa i {\ps_o} {\ps_o}) } - \xcancel{ \wh E_i (\GGa {\ps_o} {\ps_o} {\ps_o}) } - \xcancel{ \GGa {\ps_o} {\ps_o}  \ell \GGa i \ell  {\ps_o}}  - \xcancel{ \GGa {\ps_o} {\ps_o}  {\ps_o} \GGa i {\ps_o}  {\ps_o} } - \xcancel{ \GGa {\ps_o} {\ps_o}  {\qs_o} \GGa i  {\qs_o}  {\ps_o}}  + \\
\nonumber &\hskip 1.5 cm +  \xcancel{ \GGa i {\ps_o} \ell \GGa {\ps_o}  \ell {\ps_o}} + \xcancel{  \GGa i {\ps_o} {\ps_o} \GGa {\ps_o} {\ps_o} {\ps_o} } + \xcancel{ \GGa i {\ps_o} {\qs_o} \GGa  {\ps_o} {\qs_o} {\ps_o} } +\\
\nonumber & \hskip 1.5 cm+ \xcancel{\qs_o (\GGa i {\ps_o} {\qs_o}) } -  \xcancel{ \wh E_i (\GGa {\qs_o} {\ps_o} {\qs_o})}  -   \xcancel{\GGa {\qs_o} {\ps_o}  \ell \GGa i  \ell  {\qs_o}}  - \xcancel{ \GGa {\qs_o} {\ps_o}  {\ps_o} \GGa i {\ps_o}  {\qs_o} } - \xcancel{ \GGa {\qs_o} {\ps_o}  {\qs_o} \GGa i {\qs_o}  {\qs_o}} +  \\
& \hskip 1.5 cm + \xcancel{ \GGa i {\ps_o} \ell \GGa  {\qs_o} \ell {\qs_o}}  +  \xcancel{ \GGa i {\ps_o} {\ps_o} \GGa {\qs_o}   {\ps_o} {\qs_o}}  +   \xcancel{\GGa i {\ps_o} {\qs_o} \GGa {\qs_o} {\qs_o} {\qs_o} }= 0
\label{6.21}
\end{align}
\begin{align}\nonumber  & \Ric_{i{\qs_o}}\  {=} \ \RR m i {\qs_o} m + \RR {\ps_o} i {\qs_o} {\ps_o} + \RR {\qs_o} i {\qs_o} {\qs_o}  = \\
\nonumber  & \hskip 1 cm {=} \   \wh E_m (\GGa i {\qs_o} m) -  \wh E_i (\GGa m {\qs_o} m) -  \GGa m {\qs_o}  \ell \GGa i  \ell  m   - \xcancel{ \GGa m {\qs_o}  {\ps_o} \GGa  i {\ps_o}  m }  - \xcancel{ \GGa m {\qs_o}  {\qs_o} \GGa i  {\qs_o} m } +\\
\nonumber  &\hskip 1.5 cm  +   \GGa i {\qs_o} \ell  \GGa m \ell  m  +   \GGa i {\qs_o} {\ps_o}  \GGa m {\ps_o}  m +  \xcancel{ \GGa i {\qs_o} {\qs_o}  \GGa m {\qs_o}  m} - c^r_{mi} \GGa r {\qs_o} m   +  \xcancel{\o_{mi} \GGa {\qs_o} {\qs_o} m} +\\
\nonumber &\hskip 1.5 cm+  \ps_o (\GGa i {\qs_o} {\ps_o}) -  \wh E_i (\GGa {\ps_o} {\qs_o} {\ps_o}) -   \xcancel{\GGa {\ps_o} {\qs_o}  \ell \GGa i \ell  {\ps_o} } - \xcancel{ \GGa {\ps_o} {\qs_o}  {\ps_o} \GGa i {\ps_o}  {\ps_o} } -  \GGa {\ps_o} {\qs_o}  {\qs_o} \GGa i  {\qs_o}  {\ps_o}  + \\
\nonumber &\hskip 1.5 cm +  \xcancel{ \GGa i {\qs_o} \ell \GGa  {\ps_o} \ell {\ps_o}} + \xcancel{  \GGa i {\qs_o} {\ps_o} \GGa {\ps_o} {\ps_o} {\ps_o} }+ \xcancel{  \GGa i {\qs_o} {\qs_o} \GGa  {\ps_o}  {\qs_o} {\ps_o}} +
\end{align}
\begin{align}
\nonumber & \hskip 1.5 cm+ \xcancel{\qs_o (\GGa i {\qs_o} {\qs_o})} -  \wh E_i (\GGa {\qs_o} {\qs_o} {\qs_o}) -  \xcancel{ \GGa {\qs_o} {\qs_o}  \ell \GGa i \ell  {\qs_o}} - \xcancel{ \GGa {\qs_o} {\qs_o}  {\ps_o} \GGa i {\ps_o} {\qs_o} } -  \xcancel{\GGa {\qs_o} {\qs_o}  {\qs_o} \GGa i {\qs_o} {\qs_o} }+  \\
 & \hskip 1.5 cm + \xcancel{ \GGa i {\qs_o} \ell \GGa  {\qs_o} \ell {\qs_o} } +   \xcancel{ \GGa i {\qs_o} {\ps_o} \GGa  {\qs_o} {\ps_o} {\qs_o} }  +  \xcancel{ \GGa i {\qs_o} {\qs_o} \GGa {\qs_o} {\qs_o} {\qs_o}} = 0\\
 \label{6.22}
\end{align}
\begin{align}
\nonumber  & \Ric_{{\ps_o}{\qs_o}}\  {=} \ \RR m {\ps_o} {\qs_o} m +  \RR {\qs_o} {\ps_o} {\qs_o} {\qs_o}  = \\
\nonumber  & \hskip 1 cm {=} \   \xcancel{\wh E_m (\GGa {\ps_o} {\qs_o} m) } -  \ps_o (\GGa m {\qs_o} m) -  \GGa m {\qs_o}  \ell \GGa  {\ps_o} \ell m   -  \xcancel{\GGa m {\qs_o}  {\ps_o} \GGa {\ps_o} {\ps_o} m}   -  \xcancel{\GGa m {\qs_o}  {\qs_o} \GGa {\ps_o}  {\qs_o} m} +\\
\nonumber  &\hskip 1.5 cm  +  \xcancel{  \GGa {\ps_o} {\qs_o} \ell  \GGa m \ell  m}  +   \GGa {\ps_o} {\qs_o} {\ps_o}  \GGa m {\ps_o}  m +    \xcancel{\GGa {\ps_o} {\qs_o} {\qs_o}  \GGa m {\qs_o}  m }  +\\
\nonumber & \hskip 1.5 cm+  \xcancel{ \qs_o (\GGa {\ps_o} {\qs_o} {\qs_o})} -   \ps_o (\GGa {\qs_o} {\qs_o} {\qs_o}) - \xcancel{ \GGa {\qs_o} {\qs_o}  \ell \GGa  {\ps_o} \ell {\qs_o} } -  \xcancel{\GGa {\qs_o} {\qs_o}  {\ps_o} \GGa {\ps_o} {\ps_o} {\qs_o}} -   \xcancel{\GGa {\qs_o} {\qs_o}  {\qs_o} \GGa {\ps_o}  {\qs_o} {\qs_o}} +  \\
 & \hskip 1.5 cm + \xcancel{ \GGa {\ps_o} {\qs_o} \ell \GGa  {\qs_o} \ell {\qs_o}}  +    \xcancel{\GGa {\ps_o} {\qs_o} {\ps_o} \GGa  {\qs_o} {\ps_o} {\qs_o}}  +   \xcancel{ \GGa {\ps_o} {\qs_o} {\qs_o} \GGa {\qs_o} {\qs_o} {\qs_o}} = \frac{\Lambda}{2}
  \label{6.23}
  \end{align}
  \begin{align}
   \nonumber & \Ric_{{\ps_o}\ps_o}   {=} \ \RR m {\ps_o} {\ps_o} m +  \RR {\qs_o} {\ps_o} {\ps_o} {\qs_o}  = \\
\nonumber  & \hskip 1 cm {=} \  \xcancel{ \wh E_m (\GGa {\ps_o} {\ps_o} m) } -  \ps_o (\GGa m {\ps_o} m) -  \GGa m {\ps_o}  \ell \GGa {\ps_o} \ell m   - \xcancel{ \GGa m {\ps_o}  {\ps_o} \GGa {\ps_o} {\ps_o} m }  - \xcancel{ \GGa m {\ps_o}  {\qs_o} \GGa  {\ps_o} {\qs_o} m} +\\
\nonumber  &\hskip 1.5 cm  +  \xcancel{ \GGa {\ps_o} {\ps_o} \ell  \GGa m  \ell  m } +  \xcancel{ \GGa {\ps_o} {\ps_o} {\ps_o}  \GGa m {\ps_o}  m} + \xcancel{  \GGa {\ps_o} {\ps_o} {\qs_o}  \GGa  m {\qs_o}  m } +\\
\nonumber & \hskip 1.5 cm+ \xcancel{\qs_o (\GGa {\ps_o} {\ps_o} {\qs_o})} -   \xcancel{ \ps_o (\GGa {\qs_o} {\ps_o} {\qs_o})} -  \xcancel{ \GGa {\qs_o} {\ps_o}  \ell \GGa  {\ps_o} \ell {\qs_o}}  - \xcancel{ \GGa {\qs_o} {\ps_o}  {\ps_o} \GGa {\ps_o} {\ps_o} {\qs_o}} -  \xcancel{ \GGa {\qs_o} {\ps_o}  {\qs_o} \GGa  {\ps_o}  {\qs_o} {\qs_o} }+  \\
 & \hskip 1.5 cm + \xcancel{ \GGa {\ps_o} {\ps_o} \ell \GGa  {\qs_o} \ell  {\qs_o} } +  \xcancel{ \GGa {\ps_o} {\ps_o} {\ps_o} \GGa {\qs_o}  {\ps_o} {\qs_o} } +  \xcancel{ \GGa {\ps_o} {\ps_o} {\qs_o} \GGa {\qs_o} {\qs_o} {\qs_o} }= 0 \label{6.24}
\end{align}
\begin{align}\nonumber  & \Ric_{{\qs_o}{\qs_o}}\  {=} \ \RR m {\qs_o} {\qs_o} m + \RR {\ps_o} {\qs_o} {\qs_o} {\ps_o}  = \\
\nonumber  & \hskip 1 cm {=} \   \wh E_m (\GGa {\qs_o} {\qs_o} m) -  \qs_o (\GGa m {\qs_o} m) -  \GGa m {\qs_o}  \ell \GGa  {\qs_o} \ell m   -  \xcancel{ \GGa m {\qs_o}  {\ps_o} \GGa {\qs_o} {\ps_o}  m}   -  \xcancel{\GGa m {\qs_o}  {\qs_o} \GGa {\qs_o} {\qs_o} m} +\\
\nonumber  &\hskip 1.5 cm  +   \GGa {\qs_o} {\qs_o} \ell  \GGa  m \ell  m  +   \GGa {\qs_o} {\qs_o} {\ps_o}  \GGa m  {\ps_o} m +   \GGa {\qs_o} {\qs_o} {\qs_o}  \GGa m {\qs_o}  m  +\\
\nonumber &\hskip 1.5 cm+  \ps_o (\GGa {\qs_o} {\qs_o} {\ps_o}) -  \qs_o (\GGa {\ps_o} {\qs_o} {\ps_o}) -  \xcancel{ \GGa {\ps_o} {\qs_o}  \ell \GGa  {\qs_o} \ell {\ps_o} }  -  \GGa {\ps_o} {\qs_o}  {\ps_o} \GGa  {\qs_o} {\ps_o} {\ps_o}  -   \xcancel{\GGa {\ps_o} {\qs_o}  {\qs_o} \GGa {\qs_o} {\qs_o} {\ps_o} } + \\
 &\hskip 1.5 cm +   \xcancel{\GGa {\qs_o} {\qs_o} \ell \GGa  {\ps_o} \ell {\ps_o}} +  \xcancel{\GGa {\qs_o} {\qs_o} {\ps_o} \GGa {\ps_o} {\ps_o} {\ps_o} }+  \GGa {\qs_o} {\qs_o} {\qs_o} \GGa {\ps_o} {\qs_o}  {\ps_o} = \Lambda \wt \b
 \label{6.25}
\end{align}
We  decompose   this   system  into   three subsystems, namely:
\begin{itemize}[leftmargin = 25 pt]
\item[(a)] the set of equations  \eqref{6.23} -- \eqref{6.25},    concerning  the  Ricci curvatures  $\Ric_{\ps_o  \ps_o}$, $\Ric_{\ps_o \qs_o} = \Ric_{\qs_o  \ps_o}$ and $\Ric_{\qs_o \qs_o}$;
\item[(b)] the   equations   \eqref{6.21} and \eqref{6.22}, concerning the ``mixed''   curvatures $ \Ric_{{\ps_o}i}$ and $\Ric_{{\qs_o} j}$;
\item[(c)]  the equations \eqref{6.20}.
\end{itemize}
 These three subsystems are in general tightly coupled.  However under the assumption
  \beq \label{ansatz} \wh E_i(\s) = \wh E_i(\wt \b) = 0 \ ,\qquad \text{for each}\ \ 1 \leq i \leq n-2\ .\eeq
  the system becomes much   more treatable and can be solved.
  Note that these equations  together with  \eqref{62}  imply 
 \beq \label{ansatz-bis} \qs_o(\s) = \qs_o(\wt \b) = 0 \eeq
 and hence that   the functions $\s, \wt \b: M = S \times \bR \to \bR$   depend only on the coordinate of the fiber $\bR$.
 If we now assume  \eqref{ansatz} (and its consequence \eqref{ansatz-bis}),  we have that:
\begin{itemize}[leftmargin = 18pt]
\item[(1)] The Christoffel symbols $\GGa i j m$   coincide with the functions   $\Ga i j k =  g^{mk} g_o(\n^o_{E_i}   E_j, E_k)$ and they are equal to   the (lifts to $M$ of the)  Christoffel symbols of the Levi-Civita connection $\nabla^o$ of  the K\"ahler  manifold $(N, g_o, J)$ with respect to the frame field $(E_i)$. In particular, the functions $\wh{\Ric_{ij}}$ defined in \eqref{hatricci} coincide with the (lifts to $M$ of the) components  of the Ricci curvature $\Ric^N$ of  $(N,  J, g_o)$.
\item[(2)] The equations \eqref{6.21} and \eqref{6.22} reduce to
  \begin{align}
\nonumber   &   \frac{1}{4 \s}  \bigg(- \wh E_m ( J^m_i )  + \wh E_i(J^m_m)  + \Ga i m r  J^m_r  - \Ga m  r m J^r_i  + c_{m i}^\ell J_\ell^m \bigg)
 -\frac{\ps_o(\s)}{2 \s} \bigg(  \Ga i m m - \Ga m i m + c_{m i}^m\bigg)   = \\
 & \hskip 0.5 cm =  \frac{1}{4 \s}  \bigg(  -  (\n_{E_m} J)^m_i  + (\n_{E_i}J)^m_m \bigg)  + \left(   \frac{J_\ell^m }{4 \s}   -\frac{\ps_o(\s) \d_\ell^m}{2 \s} \right) \big(\Ga i m \ell  -\Ga m i \ell + c_{m i}^\ell\big) =  0
\label{6.21-bis}
\end{align}
\begin{align}
\nonumber   &  \frac{\wt \b }{2 \s}  \bigg( -\wh E_m ( J^m_i )  + \wh E_i(J^m_m)  + \Ga i m r  J^m_r   - \Ga m  r m J^r_i  + c_{m i}^\ell J_\ell^m \bigg)
  = \\
 & \hskip 3 cm =   \frac{\wt \b }{2 \s}  \bigg(   - (\n^o_{E_m} J)^m_i   + (\n^o_{E_i}J)^m_m\bigg)   +  \frac{J_\ell^m }{4 \s}  \big(\Ga i m \ell  -\Ga m i \ell + c_{m i}^\ell\big) =  0
 \label{6.22-bis}
\end{align}
  \end{itemize}
 We  observe that,  being the Levi-Civita connection $\n^o$ of the K\"ahler manifold  $(N, g_o, J)$ with trivial torsion (thus, $\Ga i m \ell  -\Ga m i \ell + c_{m i}^\ell \equiv 0$) and  such that  $\n^o J = 0$, it follows that  \eqref{6.21-bis} and \eqref{6.22-bis}  are identically satisfied. This means  that the  query for  Einstein metrics  now reduces just to  looking for solutions
 to  the  equations   \eqref{6.20} and  \eqref{6.23} -- \eqref{6.25}.
\par
\medskip
Since  $J$ is a complex structure, we have   $J^m_m = \operatorname{tr}(J) = 0$ and $J^m_\ell J^\ell_m  =\operatorname{tr}(J^2) = - (n-2)$. Using this  together   with the identity $\wh \Ric_{ij} = \Ric^N_{ij}$ and the  expressions \eqref{6.12} -- \eqref{6.17},  we  see  that  \eqref{6.20} is  equivalent to
\beq
\Ric^N_{ij} = -   \left(2\ps_o(\wt \b)\ps_o(\s)+\left( \frac{n-4}{2\s}(\ps_o(\s))^2 -\frac{1}{4\s}+ \ps_o(\ps_o(\s))\right)2\wt \b -  \s \Lambda\right) g_{ij}
 \label{6.20-bis}
\eeq
and that  the equations  \eqref{6.23} -- \eqref{6.25} are equivalent to
  \begin{align}
  &2  \ps_o(\ps_o(\wt \b))+\frac{(n-2)}{\s}\ps_o(\wt \b)\ps_o(\s)+\frac{(n-2)}{4\s^2}\wt \b=\Lambda \ ,
  \label{6.23-bis}  \\[15 pt]
&  \frac{ n-2}{4 \s^2} \bigg( -2 \s  \ps_o (  \ps_o(\s))+ (\ps_o(\s))^2  +  \frac{1}{4}   \bigg)  =  0 \ ,\label{6.24-bis} \\[15pt]
 & 2 \wt \b \left ( 2\ps_o(\ps_o(\wt \b))+\frac{(n-2)}{\s}\ps_o(\wt \b)\ps_o(\s)+\frac{(n-2)}{4\s^2}\wt \b\right)= 2 \Lambda \wt \b \ .
 \label{6.25-bis}
\end{align}
We observe that, since the left hand side  in  \eqref{6.20-bis} is independent of the fiber coordinates of the bundle $\pi^S \circ \pi: M \to N$, such equation might be satisfied only if
$(N, g_o, J)$ is K\"ahler-Einstein, i.e. $\Ric = \L^o g$ for some  Einstein constant $\L^o$, and $\s$ and $\wt \b$ satisfy the equation
\beq \label{6.20-ter} 2  \ps_o(\wt \b)\ps_o(\s)+\left( \frac{n-4}{2\s}(\ps_o(\s))^2 -\frac{1}{4\s}+ \ps_o(\ps_o(\s))\right)2 \wt \b - \s \Lambda + \Lambda^o = 0\ .\eeq
Moreover, since \eqref{6.25-bis} is manifestly implied by \eqref{6.23-bis} and we  look for solutions in which  $\s > 0$ at all points,  it now suffices to find  solutions $\s > 0$ and $\wt \b$ of the system of the just  {\it three} ordinary differential equations \eqref{6.20-ter}, \eqref{6.23-bis} and \eqref{6.24-bis}. \par
\smallskip
Now we  show that  the equation \eqref{6.23-bis}  follows from the others.   Indeed, using the fact that $\s >0$ and differentiating  the whole term inside parentheses in \eqref{6.24-bis} along $\ps_o$, we immediately get that  $\s$ satisfies
\beq\label{6.26} \ps_o( \ps_o (  \ps_o(\s))) = 0\ . \eeq
 Using this property and  replacing  $\ps_o (  \ps_o(\s))$  by the expression $\ps_o (  \ps_o(\s)) =   \frac{1}{2 \s} \big((\ps_o(\s))^2  +  \frac{1}{4}\big)$ which follows
 from   \eqref{6.24-bis}, one can check that the  derivative along $\ps_o$ of  the equation \eqref{6.20-ter}   is the multiple by $\ps_o(\s)$ of the equation \eqref{6.23-bis}. This proves  that the latter is  implied  by  \eqref{6.20-ter} in case $\ps_o(\s) \neq 0$.
 \par
\smallskip
This means that we are now    reduced to find  solutions $\s > 0$ and $\wt \b$  to the o.d.e. system
\begin{align}
\label{prima} &2  \ps_o(\wt \b)\ps_o(\s)+\left( \frac{n-4}{2\s}(\ps_o(\s))^2 -\frac{1}{4\s}+ \ps_o(\ps_o(\s))\right) 2 \wt \b - \s \L + \L^o = 0\ ,\\
 \label{seconda} & -2 \s  \ps_o (  \ps_o(\s))+ (\ps_o(\s))^2  +  \frac{1}{4}  = 0\ .
\end{align}
 We now observe that  \eqref{6.26} shows that each solution  $\s = \s(t)$ to \eqref{seconda} is  a quadratic polynomial  $ \s(t) = C_2 t^2 + C_1 t + C_0$.  Actually, it is equivalent to say that  such positive quadratic polynomial  has 
 discriminant equal to $- \frac{1}{4}$.  So, by an appropriate coordinate change $t \mapsto t + c$,  we may always assume that  it has the simpler  form
\beq
  \label{solutionsigma1} \s(t) =  \frac{1}{16 C}  t^2 + C \qquad \text{for some}\ C > 0\ .\eeq
  From this we get that  \eqref{prima} is equivalent  to the first order linear differential equation
\beq \label{eqbeta}
\frac{t}{4 C} \frac{d \wt \b}{dt}  {+}\left( \frac{n-4}{\frac{t^2}{16 C} + C}\left(\frac{t}{8 C} \right)^2 -\frac{1}{2\left(\frac{t^2}{16 C} + C\right)}+  \frac{1}{4C}\right) \wt \b {-} \left(\frac{t^2}{16 C} + C\right) \L {+} \L^o = 0\ .
\eeq
The general solutions on $(0, + \infty)$ of this equation have the form
 \begin{multline} \label{eqtheorem-1}
 \wt \b(t) \=  e^{-\int_1^t a(s) d s} \left(B_0 -  \int_1^t b(s) e^{\int_1^s a(\check s) d \check s} d s\right)\qquad \text{with}\\
  a (t) \= \frac{8 C}{t}\left( \frac{n-4}{2\left(\frac{t^2}{16 C} + C\right)}\left(\frac{t}{8 C} \right)^2 -\frac{1}{4\left(\frac{t^2}{16 C} + C\right)}+  \frac{1}{8C}\right)  =\frac{(n-3)t^2-16C^2}{t^3+16C^2t}\ \text{and}\ \\
 b(t) \= - \frac{4 C}{t} \left(\frac{t^2}{16 C} + C\right) \L + \frac{4 C}{t}\L_o\ . \end{multline}
We now observe that  $e^{\int_1^t a(s)ds} = B_1 \frac{(t^2+16C^2)^{\frac{n}{2}-1}}{t}$  for some  $B_1 \in \bR$ and that $b(s)e^{\int_1^s a(\check s) d\check s} =B_1 \frac{\left(16 C^2+s^2\right)^{\frac{n}{2}-1}  \left(16 C \Lambda _0- \Lambda \left(16 C^2+s^2\right)\right)}{4 s^2} $. 
Hence,  since $n$ is even, it is a rational function of $s^2$. Plugging these expressions into \eqref{eqtheorem-1}, we  see that $\wt \b(t)$ is a rational function that  is well defined  at $t= 0$ and hence on the whole real axis. Setting $B = B_0/B_1$,  we get \eqref{eqtheorem}. 
\end{pf}

\begin{rem}    As a corollary of the   proof,  we have  that if  the  quantisable  K\"ahler manifold $(N, J,  g^o)$ is  not Einstein,   there is no  firmly compatible  Einstein metric on the  trivial $\bR$-bundle $M = S \times \bR$.
\end{rem}
\par
As we will shortly see,   the classical Taub-NUT metrics can be considered as   firmly compatible metrics \eqref{eqtheorem} on the $4$-dimensional RT bundle $\pi: M = S^3 \times \bR \to \bR$, where $S^3$ is considered as a Sasaki  manifold over the round sphere $S^2 = \bC P^1$. Due to this,   we call  the metrics of Theorem \ref{main}   {\it of Taub-NUT type}. \par
\medskip
\subsection{$4$-dimensional Taub-NUT metrics and  higher dimensional analogues}
Let $(S, \cD = \ker \theta, J)$ and $(N, J, g_o)$ be a regular Sasaki manifold and its associated K\"ahler-Einstein manifold, respectively,  as in Theorem \ref{main}.
Consider  a (local) trivialisation $S|_{\cU} = \cU \times A$  on  some open set $\cU \subset N$ and denote by  $u$ a coordinate  for the  fiber $A = \bR, S^1$ of the Sasaki manifold.
In this way,  the contact form $\theta$ takes the form  $\theta = du + \h$ for some appropriate $1$-form $\h$ in $\O^1(\cU)$.   Then, using  the   coordinates $(t,u)$  for the  fiber $\bR \times A $ of  the (trivialised) bundle  $\pi^S \circ \pi: M|_\cU =  S|_\cU \times \bR \to \cU \subset N$,  the   metrics of Theorem \ref{main}  read  as
\begin{multline*} g =   \s(t)  g_o +  \wt \b(t) \left\{ \big(du + \h\big) \vee \bigg( \frac{1}{ \wt \b(t)} dt +   du + \h\bigg) \right\}= \\
=  \s(t)  g_o +  \wt \b(t) \left( \bigg(du  +  \h +   \frac{1}{2\wt \b(t)} dt\bigg) -   \frac{1}{2\wt \b(t)} dt \right)   \vee \left(   \bigg(du  +  \h + \frac{1}{2\wt \b(t)} dt \bigg) +  \frac{1}{2\wt \b(t)} dt   \right)
\end{multline*}
\begin{multline}\label{coorTNUT}
=  \s(t)  g_o + \wt \b(t)  \bigg(du  +  \h +   \frac{1}{2\wt \b(t)} dt\bigg)^2  -   \frac{1}{4\wt \b(t)} dt ^2 = \\
= \s(t)  g_o + \wt \b(t) \bigg(\frac{1}{\L_o} d v +  \h\bigg)^2  -   \frac{1}{4\wt \b(t)} dt ^2
\end{multline}
where, in the last step,  we  replaced the fiber coordinates $(t, u)$  by
\beq (t,u) \longmapsto \left(t, v \= \L_o\bigg(u + \int_0^t \frac{1}{2 \wt \b(s)} ds\bigg)\right)\ .\eeq
 If we now take   $N = \bC P^1 = S^2$  as K\"ahler manifold and  consider its standard spherical coordinates $(\phi, \psi)$,   the
 round metric $g_o$ of  constant curvature $\k = \L_o$  on $N = S^2$ and the associated  contact form $\theta$ on the associated Sasaki manifold $S^3$ (which  fibers on $S^2$ by means of the Hopf map) have the   coordinate expressions
\beq \label{5.39}  g_o = \frac{1}{\L_o} (d \psi ^2+ \sin \psi^2 d \phi^2)\ ,\qquad \theta = du + \h = d u + \frac{1}{\L_o} \cos \psi d \phi\ .\eeq
On the other hand the functions \eqref{eqtheorem}   corresponding to   Ricci flat metrics  ($\L = 0$)  are
 \beq  \label{5.40}  \wt \b(t) = \L_o\left( \check B  \frac{t}{t^2 + 16 C^2}  -  4 C \frac{ t^2 - 16 C^2}{t^2 + 16 C^2}\right)\qquad \text{with}\ \  \check B := \frac{ B}{2\L_o} (1 + 16 C^2) + 4 C  (1 - 16 C^2)\ .\eeq
Plugging \eqref{5.39} and \eqref{5.40} into \eqref{coorTNUT},  we obtain  the following coordinate expression for the Ricci flat metrics  on the $4$-manifold $M = S^3 \times \bR = \bR^4 \setminus\{0\}$:
\begin{multline} \label{coorTNUT-1}
g
= \frac{1}{\L_o} \left(\frac{1}{16 C}  t^2 + C\right)  (d \psi ^2+ \sin \psi^2 d \phi^2) +\\
+   \frac{1}{\L_o}  \frac{ \check Bt - 4 C  t^2 + 64 C^3}{t^2 + 16 C^2}  \left(d v +   \cos \psi d \phi\right)^2
-  \frac{1}{\L_o}   \frac{t^2 + 16 C^2}{4 \left(\check B t   -  4 C t^2 + 64 C^3)\right)} dt ^2
\end{multline}
If we now set
\begin{multline}\ell \= \sqrt{C} \ .\qquad m\=  \frac{\check B}{32 \ell^3} =  \frac{1}{32 C^{\frac{3}{2}}} \left( \frac{ B}{2\L_o} (1 + 16 C^2) + 4 C  (1 - 16 C^2)\right)\\
\text{and apply the coordinate  change}\ \
t \longmapsto \check t = \frac{t}{4 \ell} = \frac{t}{4 \sqrt{C}}\end{multline}
 the metrics \eqref{coorTNUT-1} take the very familiar coordinate expression of the (rescaled)  Taub-NUT metrics
 \begin{multline}   g\= \frac{1}{\L_o} \bigg\{ (\check t^2 + \ell^2)(d\psi ^{2}+\sin^2 \psi d\phi ^{2}) +\\
+     \frac{ 2 m  \check t  + \ell^2 -   \check t^2}{   \check t^2 +  \ell^2} 4 \ell^2 \left(d v +   \cos \psi d \phi\right)^2  -    \frac{ \check t^2 +  \ell^2}{  2 m \check t   +  \ell^2 -  \check  t^2} d\check t ^2\bigg\} \ . \end{multline}
Many other  $4$-dimensional  Einstein metrics   can be determined in exactly  the same way: it suffices to   impose a different value  $\L \neq 0$ for the desired Einstein constant  and/or
to replace the K\"ahler manifold $N = S^2$ by some other  compact Riemann surface, as e.g.    $T^2 = S^1 \times S^1$ or  a  compact quotient of the unit disk $\D \subset \bC$,  equipped with    some metric of constant curvature. \par
\smallskip
In order to  generate explicit examples of higher dimensional  Einstein metrics of Taub-NUT type,   one might follow the same  procedure as above,  starting, for instance,  from some higher dimensional  homogeneous   flag manifolds $N = G^\bC/P$ of a complex  semisimple Lie group,  equipped with its unique (up to a homothety) invariant K\"ahler-Einstein metric of  positive Einstein constant $\L_o> 0$.    For instance, one might consider the $4$-dimensional  manifolds $N = \bC P^2$ or $N =\bC P^1 \times \bC P^1$ equipped   with the  Fubini-Study metric or the cartesian product of two round metrics, respectively. In all these  cases,  for any choice of a constant $\L \in \bR$ and of a Sasaki manifold $(S, \cD = \theta, J)$ projecting onto $N$  (about the regular Sasaki manifolds that are associated with  the  compact flag manifolds, see for instance \cite{ACHK} and references therein),  Theorem \ref{main} and the
above discussion   yield  explicit coordinate expressions for   Lorentzian Einstein metrics with  any prescribed Einstein constant $\L$   over  $M = S \times \bR$.\par
\par \medskip
\appendix
\section{The Christoffel symbols of the Levi-Civita connection of  a Lorentzian  metric  of  K\"ahler-Sasaki type}  \label{Tedious}
In this appendix, we give the explicit expressions of the Christoffel symbols of the Levi-Civita connection of a compatible metric on the trivial RT bundle $\pi: M = S \times \bR \to S$
as described in \S \ref{severe}. We compute them in terms of the frame field   $(\ps_o, \wh E_i,  \qs_o)$ and its dual coframe field $(\ps_o^*, \wh E^i, \qs^*_o)$, associated with a local frame field $(E_i)$. More precisely, we  list here the real  functions $\GGa A B C$ which determine the   covariant derivatives
$\n_{X_A} X_ B = \GGa A B C X_C$ for any choice of  two vector fields  $X_A, X_B$ of the tuple  $(\ps_o, \wh E_i,  \qs_o)$ for a metric $g$ of the form (see \eqref{buona})
\beq  \label{buona-ter} g  = \s\left( g_{ij} \wh E^i\vee \wh E^j + \qs^*_o \vee \left(\a\ps^*_o   + \g^i g_{ij} \wh E^j  + \b  \qs_o^*\right) \right) \ ,\qquad g_{ij} = g_o(E_i, E_j)\ .
\eeq
First we  compute the Christoffel symbols under the simplifying assumption
 $\s \equiv 1$. Then   one can directly check  that the elements of the dual coframe field
 $(\ps_o^*,\wh E^1, \ldots, \wh E^{n-2},  \qs_o^*)$ satisfy the identities
 \begin{multline}\label{dualcoframe}  \wh E^i = g\bigg(g^{ik} \wh E_k - \frac{\g^i}{\a} \ps_o, \cdot\bigg) \ ,\quad
 \ps_o^* = g\bigg(\frac{2}{\a} \qs_o +  \frac{1}{\a^2} \left( \g^m \g^k  g_{mk} -  4 \b \right) \ps_o - \frac{ \g^m}{\a} \wh E_m, \cdot\bigg)\ ,\\
  \qs_o^* = g\bigg(\frac{2}{\a} \ps_o, \cdot\bigg)\ .
\end{multline}
Using  these  and the expansion $X = \wh E^i(X) \wh E_i + \ps_o^*(X) \ps_o + \qs_o^*(X) \qs_o$ of any vector field $X \in \gX(M)$ in terms  of our  frame field,  we may   compute the  covariant derivatives $\n_{X_A} X_B$  for any  pair $X_A, X_B$ of vector fields of  the frame  $\{\wh E_i,  \ps_o, \qs_o\}$ using
  Koszul's formula
\begin{multline} \label{Koszul}  g(\n_X Y, Z) = \frac{1}{2} \bigg( X(g(Y, Z)) + Y(g( X, Z)) - Z(g(X, Y)) - \\
  - g([X, Z], Y) - g(  [Y, Z], X) + g([X, Y],Z)\bigg) \ .\end{multline}
  Note also that for  any pair of commuting  vector  fields  $X_A$, $X_B$ one has $\n_{X_A} X_B = \n_{X_B} X_A$  because   $\n$ has  trivial torsion.
From   this information and  by  some (long and tedious, but very   straightforward)  computations we get the following  list of  covariant derivatives in the case $\s \equiv 1$. In the next formulas we use the notation  $g_{ij} \= g_o(E_i, E_j)$,    $\o_{ij} \= g_o(E_i, J E_j)$, $[E_i, E_j] = c_{ij}^k E_k$  and we set
  \begin{align*}
 &   \ S_{ij|k} \=    \frac{1}{2} \bigg(  \frac{\g^\ell}{2}g_o(E_i, JE_k) g_o(E_\ell, E_j)  + \frac{\g^\ell}{2} g_o(E_j, J E_k)  g_o(E_\ell, E_i)+\\
& \hskip 9 cm -   \frac{\g^\ell}{2} g_o(E_i, JE_j)  g_o(E_\ell,E_k)\bigg) \ .
\end{align*}
Here is the list:
 \begin{multline} \label{731}
\n_{\wh E_i} \wh E_j {=}  \left(
  g^{mk} g_o(\n^o_{E_i}   E_j, E_k)  +  g^{mk}S_{ij|k}
  +  \frac{\g^m \o_{ij} }{4} \right) \wh E_m
+  \left(    \frac{1}{2 \a} \wh E_i(\g^k g_{jk} )+   \frac{1}{2 \a}\wh E_j(\g^k g_{ik} )  - \right.\\
\left.  -   \frac{1}{4 \a} \g^m \g^k  g_{mk} \o_{ij}   - \frac{ \g^m}{\a} g_o(\n^o_{E_i}   E_j, E_m)  - \frac{ \g^m}{\a}  S_{ij|m} \right)\ps_o -   \frac{  \o_{ij}}{2} \qs_o\ .
 \end{multline}
\beq
 \n_{\wh E_i} \ps_o =   \frac{\a g^{mk} \o_{ik} }{4}  \wh E_m
+  \left(    \frac{1}{2\a}  \wh E_i(\a) +  \frac{1}{2\a}  \ps_o( \g^k) g_{ik} - \frac{ \g^m \o_{im}}{4}
 \right)\ps_o  \hskip 2 cm
\eeq
 \begin{multline}
 \n_{\wh E_i} \qs_o   = \Bigg(\frac{g^{mk}}{4} \wh E_i (\g^t g_{tk})- \frac{g^{mk}}{4}\wh E_k (\g^t g_{ti}) + \frac{g^{mk}}{2}\b \o_{ik} - \frac{\g^\ell}{4} c_{i r}^t g_{t \ell} g^{m r} -\\
 \hskip 7 cm - \frac{\g^m}{4\a}\wh E_i(\a)  +  \frac{\g^m}{4\a}\ps_o(\g^t)g_{it}\Bigg)\wh E_m+\\
 +\bigg( \frac{1}{\a}\wh E_i(\b)+ \frac{1}{4\a^2} \g^m \g^k  g_{mk}  \wh E_i(\a)  -  \frac{1}{4\a^2} \g^m \g^k  g_{mk} \ps_o(\g^t)g_{it}  -  \frac{\b}{\a^2}\wh E_i(\a)  + \frac{1}{\a^2}\b\ps_o(\g^t)g_{it} -\\
  -\frac{\g^m}{4\a}\wh E_i(\g^tg_{tm}) + \frac{\g^m}{4\a}\wh E_m(\g^tg_{it})  - \frac{\g^m}{2\a}\b\o_{im}\bigg)\ps_o +\Bigg( \frac{1}{2\a}\wh E_i(\a) - \frac{1}{2\a}\ps_o(\g^t)g_{it}\Bigg)\qs_o
\end{multline}
\beq
  \n_{\ps_o} \wh E_i =  \frac{\a g^{mk}}{4}\o_{ik}\wh E_m +\Bigg( \frac{1}{2\a}\wh E_i(\a)+ \frac{1}{2\a}\ps_o(\g^t)g_{it} - \frac{\g^m}{4}\o_{im} \Bigg)\ps_o 
 \eeq
\beq
  \n_{\ps_o} \ps_o =
 \left(\ps_o(\log \a) \right)\ps_o
\eeq
\begin{multline}
  \n_{\ps_o} \qs_o =\bigg(  \frac{1}{4}\ps_o(\g^m) -  \frac{g^{mk}}{4}\wh E_k(\a)\Bigg)  \wh E_m+\Bigg( \frac{1}{\a}\ps_o(\b) -\frac{\g^m}{4\a}\ps_o(\g^i)g_{im}+ \frac{\g^m}{4\a}\wh E_m(\a)   \bigg)\ps_o
\end{multline}
\begin{multline}
 \n_{\qs_o} \wh E_i = \Bigg(\frac{g^{mk}}{4} \wh E_i (\g^t g_{tk})- \frac{g^{mk}}{4}\wh E_k (\g^t g_{ti}) + \frac{g^{mk}}{2}\b \o_{ik}  - \frac{\g^\ell}{4} c_{i r}^t g_{t \ell} g^{m r} -\\
 \hskip 7 cm - \frac{\g^m}{4\a}\wh E_i(\a)  +  \frac{\g^m}{4\a}\ps_o(\g^t)g_{it}\Bigg)\wh E_m+\\
 +\bigg( \frac{1}{\a}\wh E_i(\b)+ \frac{1}{4\a^2} \g^m \g^k  g_{mk}  \wh E_i(\a)  -  \frac{1}{4\a^2} \g^m \g^k  g_{mk} \ps_o(\g^t)g_{it}  -  \frac{\b}{\a^2}\wh E_i(\a)  + \frac{1}{\a^2}\b\ps_o(\g^t)g_{it} -\\
  -\frac{\g^m}{4\a}\wh E_i(\g^tg_{tm}) + \frac{\g^m}{4\a}\wh E_m(\g^tg_{it})  - \frac{\g^m}{2\a}\b\o_{im}\bigg)\ps_o +\Bigg( \frac{1}{2\a}\wh E_i(\a) - \frac{1}{2\a}\ps_o(\g^t)g_{it}\Bigg)\qs_o
\end{multline}
\begin{multline}
  \n_{\qs_o} \ps_o =\Bigg(  \frac{g^{mk}}{4} \ps_o(\g^i)g_{ik}-  \frac{g^{mk}}{4}\wh E_k(\a)\Bigg)\wh E_m+\Bigg( \frac{1}{\a}\ps_o(\b) -\frac{\g^m}{4\a}\ps_o(\g^t)g_{tm}+\frac{\g^m}{4\a}\wh E_m(\a) \Bigg)\ps_o
\end{multline}
\ \\[-20 pt]
\begin{multline*}
  \n_{\qs_o} \qs_o = \Bigg(\frac{g^{mk}}{2}\qs_o(\g^i)g_{ik} - \frac{g^{mk}}{2}\wh E_k(\b)  -\frac{\g^m}{2\a} \qs_o(\a) +\frac{\g^m}{2\a}\ps_o(\b)\Bigg)\wh E_m+\\
  +\Bigg(\frac{1}{ \a} \qs_o(\b)+   \frac{1}{2\a^2}\g^m \g^k  g_{mk} \qs_o(\a)-  \frac{1}{2\a^2}\g^m \g^k  g_{mk}\ps_o(\b)   -  \frac{2}{\a^2}\b\qs_o(\a) +   \frac{2}{\a^2}\b\ps_o(\b) -
  \end{multline*}
  \beq  \label{739}
  -\frac{\g^m}{2\a} \qs_o(\g^i)g_{im} + \frac{\g^m}{2\a}\wh E_m(\b) \Bigg)\ps_o+ \left(  \frac{1}{\a}\qs_o(\a)-  \frac{1}{\a}\ps_o(\b)\right)\qs_o
\eeq
 We now  denote  by  $\wt g$ a compatible metric as above, determined by the conformal factor $\s = 1$ and  by  $g = \s \tilde{g}$ another compatible metric, which is determined by  an arbitrary conformal factor $\s> 0$.  The Levi-Civita connection $D$ of $g$ is related with the Levi-Civita connection $\n$ of $\wt g$ by  the  formula (see e.g. \cite{Be}*{Th. 1.159})
\beq  \label{confD} D_X Y = \nabla_X Y + X(\f) Y + Y(\f) X -  g(X, Y)  \grad(\f) \ ,\qquad \f \= \frac{1}{2} \log \s\ .\eeq
On the other hand,  by \eqref{dualcoframe},
\beq \grad \f = (\grad \f)^{\wh E_i} \wh E_i +  (\grad \f)^{\ps_o} \ps_o + (\grad \f)^{\qs_o} \qs_o\ ,\eeq
where
\begin{align} \label{742}
& (\grad \f)^{\wh E_i} \= g^{ik} \wh E_k(\f) - \frac{\g^i}{\a} \ps_o(\f)\ ,\\
& (\grad \f)^{\ps_o} \=  \frac{2}{\a} \qs_o(\f) +  \frac{1}{\a^2} \left(  \g^m \g^k  g_{mk} - 4\b \right) \ps_o(\f) - \frac{ \g^m}{\a} \wh E_m(\f)\ ,\\
\label{744} &  (\grad \f)^{\qs_o} \=  \frac{2}{\a} \ps_o(\f) .
\end{align}
Combining  \eqref{731} -- \eqref{739}  with \eqref{confD},  \eqref{742} -- \eqref{744}, we  see that  the Christoffel symbols   $\GGa A B C$ of  an arbitrary  compatible metric $g$ are 
\begin{align}
\label{746}
\GGa i j m &=
  g^{mk} g_o(\n^o_{E_i}   E_j, E_k)  +  g^{mk}S_{ij|k}
  +  \frac{\g^m \o_{ij} }{4}   + \frac{1}{2 \s}  \wh E_i(\s) \d_j^m +  \frac{1}{2 \s}  \wh E_j(\s) \d_i^m\nonumber \\
 &\hskip 1 cm- \frac{g_{ij}}{2 \s}   \left(  g^{mk} \wh E_k( \s) - \frac{\g^m}{\a} \ps_o(\s)\right)\\
 \GGa i j {\ps_o}  &=  \frac{1}{2 \a} \wh E_i(\g^k g_{jk} )+  \frac{1}{2 \a}\wh E_j(\g^k g_{ik} )
-  \frac{1}{4 \a} \g^m \g^k  g_{mk} \o_{ij}  - \frac{ \g^m}{\a} g_o(\n^o_{E_i}   E_j, E_m)  - \frac{ \g^m}{\a}  S_{ij|m} \nonumber\\
 &\hskip 1 cm- \frac{g_{ij}}{2 \s}  \left(\frac{2}{\a} \qs_o(\s) +  \frac{1}{\a^2} \left(  \g^m \g^k  g_{mk} - 4\b \right) \ps_o(\s) - \frac{ \g^m}{\a} \wh E_m(\s)\right)\\
 \GGa i j {\qs_o}  &= - \frac{  \o_{ij}}{2}  - \frac{ g_{ij} }{\a \s} \ps_o(\s)
  \end{align}
 \begin{align}
 \GGa i {\ps_o} m &=  \GGa {\ps_o} i m =   \frac{\a g^{mk} \o_{ik} }{4} +\frac{1}{2 \s}  \ps_o(\s) \d_i^m \\
 \GGa i {\ps_o} {\ps_o}&= \GGa {\ps_o} i {\ps_o} = \frac{1}{2\a}  \wh E_i(\a) +  \frac{1}{2\a}  \ps_o( \g^k) g_{ik} - \frac{ \g^m \o_{im}}{4} + \frac{1}{2 \s} \wh E_i(\s) \\
 \GGa i {\ps_o} {\qs_o}&=  \GGa {\ps_o} i {\qs_o} = 0\\
 \GGa i {\qs_o} m &=  \GGa  {\qs_o} i m = \frac{g^{mk} }{4}\wh E_i(\g^tg_{tk})- \frac{g^{mk} }{4}\wh E_k(\g^tg_{ti})+ \frac{g^{mk} }{2}\b\o_{ik}  - \frac{\g^\ell}{4} c_{i r}^t g_{t \ell} g^{m r} - \nonumber \\
& \hskip 0.5 cm -\frac{\g^m}{4\a}\wh E_i(\a) + \frac{\g^m}{4\a}\ps_o(\g^t )g_{ti} + \frac{1}{2 \s} \qs_o(\s)\d_i^m-   \frac{\g^t}{4 \s }  g_{ti} \left(g^{mk} \wh E_k(\s) - \frac{\g^m}{\a} \ps_o(\s)\right)
\\
 \GGa i {\qs_o} {\ps_o} &= \GGa  {\qs_o} i {\ps_o}  =   \frac{1}{\a}\wh E_i(\b)+  \frac{1}{4 \a^2} \g^m \g^k g_{mk} \wh E_i(\a) - \frac{1}{4 \a^2} \g^m \g^k g_{mk}\ps_o(\g^t )g_{it}  - \frac{1}{ \a^2} \b\wh E_i(\a) + \nonumber \\
& \hskip 1 cm +  \frac{1}{ \a^2} \b \ps_o(\g^t )g_{it} - \frac{\g^m}{4\a} \wh E_i(\g^tg_{tm}) +  \frac{\g^m}{4\a}\wh E_m(\g^tg_{it}) -  \frac{\g^m}{2\a}\b\o_{im} - \nonumber  \\
 & \hskip 1 cm-    \frac{\g^t}{4 \s} g_{ti} \left(\frac{2}{\a} \qs_o(\s) +  \frac{1}{\a^2} \left(  \g^m \g^k  g_{mk} - 4\b \right) \ps_o(\s) - \frac{ \g^m}{\a} \wh E_m(\s)\right) \\
 \GGa i {\qs_o} {\qs_o} &=  \GGa  {\qs_o} i {\qs_o} =  \frac{1}{2\a}\wh E_i(\a)- \frac{1}{2\a}\ps_o(\g^t)g_{it}  + \frac{1}{2 \s}  \wh E_i(\s)   -     \frac{\g^t  g_{ti}}{2 \a \s}  \ps_o(\s)
 \\
 \GGa {\ps_o}  {\ps_o}  m & =  0  \\ %
 \GGa {\ps_o}   {\ps_o} {\ps_o} & =  \ps_o(\log ( \a \s))\\
 \GGa {\ps_o}   {\ps_o} {\qs_o} & =  0 \\
 \GGa {\ps_o}  {\qs_o}  m & =  \GGa {\qs_o}  {\ps_o}  m =  \GGa {\qs_o}  {\ps_o}  m = \frac{1}{4}\ps_o(\g^m) - \frac{g^{mk}}{4 } \wh E_k(\a)-  \frac{\a }{4\s} \left (g^{mk} \wh E_k(\s) - \frac{\g^m}{\a} \ps_o(\s) \right)
\\
 \GGa {\ps_o}   {\qs_o} {\ps_o} & =  \GGa {\qs_o}   {\ps_o} {\ps_o} = \frac{1}{\a}\ps_o(\b) -\frac{\g^m}{4\a}\ps_o(\g^i)g_{im}+ \frac{\g^m}{4\a}\wh E_m(\a)+\frac{1}{2 \s} \qs_o(\s)-\nonumber
\\
 & \hskip 1 cm- \frac{1}{2 \s}\left(\qs_o(\s) +  \frac{1}{ 2 \a} \left(  \g^m \g^k  g_{mk} - 4\b \right) \ps_o( \s) - \frac{ \g^m}{2 } \wh E_m(\s) \right)  \\
 \GGa {\ps_o}   {\qs_o} {\qs_o} &  = \GGa {\qs_o}   {\ps_o} {\qs_o}   = 0 \\
\GGa {\qs_o}  {\qs_o}  m & = \frac{g^{mk}}{2}\qs_o(\g^i)g_{ik} - \frac{g^{mk}}{2}\wh E_k(\b)  -\frac{\g^m}{2\a} \qs_o(\a) +\frac{\g^m}{2\a}\ps_o(\b)-\nonumber\\
& \hskip 1 cm-  \frac{\b }{2 \s} \left( g^{mk} \wh E_k(\s) - \frac{\g^m}{\a} \ps_o(\s)\right) \\
\GGa {\qs_o}   {\qs_o} {\ps_o} & = \frac{1}{\a} \qs_o(\b)+   \frac{1}{2\a^2}\g^m \g^k  g_{mk} \qs_o(\a)-  \frac{1}{2\a^2}\g^m \g^k  g_{mk}\ps_o(\b)   -  \frac{2}{\a^2}\b\qs_o(\a) + \nonumber\\
& \hskip 1 cm +  \frac{2 \b}{\a^2}\ps_o(\b)
   -\frac{\g^m}{2\a} \qs_o(\g^i)g_{im} + \frac{\g^m}{2\a}\wh E_m(\b)-\nonumber \\
& \hskip 1 cm-  \frac{\b }{\s}\left(\frac{1}{\a} \qs_o(\s) +  \frac{1}{2\a^2} \left(  \g^m \g^k  g_{mk} - 4\b \right) \ps_o(\s) - \frac{ \g^m}{2 \a} \wh E_m(\s)\right)
\\
\label{772} \GGa {\qs_o}   {\qs_o} {\qs_o} & =  \frac{1}{\a}\qs_o(\a)-  \frac{1}{\a}\ps_o(\b)+\frac{1}{\s} \qs_o(\s)- \frac{\b}{ \a \s}\ps_o(\s)
  \end{align}
  \par

\vskip 1.5truecm
\hbox{\parindent=0pt\parskip=0pt
\vbox{\baselineskip 9.5 pt \hsize=3.5truein
\obeylines
{\smallsmc
Dmitri V.  Alekseevsky
Institute for Information Transmission Problems
B. Karetny per. 19
127051 Moscow
Russia
\&
University of Hradec   Kr\'alov\'e,
Faculty of Science, 
Rokitansk\'eho 62, 
500~03 Hradec Kr\'alov\'e,  
Czech Republic
%
%
%
\
}\medskip
{\smallit E-mail}\/: {\smalltt dalekseevsky@iitp.ru}
}
\vbox{\baselineskip 9.5 pt \hsize=3.1truein
\obeylines
{\smallsmc
Masoud Ganji
School of Science and Technology
University of New England,
Armidale NSW 2351
Australia
\phantom{\&}
\phantom{University of Hradec   Kr\'alov\'e,}
\phantom{Faculty of Science, Rokitansk\'eho 62, 500~03 Hradec Kr\'alov\'e,}
\phantom{Czech Republic}
\phantom{Czech Republic}
\phantom{Czech Republic}
\phantom{Czech Republic}}\medskip
{\smallit E-mail}\/: {\smalltt mganjia2@une.edu.au}
}
}
\vskip 1truecm
\hbox{\parindent=0pt\parskip=0pt
\vbox{\baselineskip 9.5 pt \hsize=3.5truein
\obeylines
{\smallsmc
Gerd Schmalz
School of Science and Technology
University of New England,
Armidale NSW 2351
Australia
\
}\medskip
{\smallit E-mail}\/: {\smalltt schmalz@une.edu.au}
}
\vbox{\baselineskip 9.5 pt \hsize=3.1truein
\obeylines
{\smallsmc
Andrea Spiro
Scuola di Scienze e Tecnologie
Universit\`a di Camerino
Via Madonna delle Carceri
I-62032 Camerino (Macerata)
Italy
}\medskip
{\smallit E-mail}\/: {\smalltt andrea.spiro@unicam.it
}
}
}
\end{document}